\RequirePackage[normalem]{ulem} 
\RequirePackage{color}\definecolor{RED}{rgb}{1,0,0}\definecolor{BLUE}{rgb}{0,0,1} 

\documentclass{elsarticle}

\journal{arXiv}

\usepackage{amsmath,stmaryrd,amsfonts,amssymb,fancyhdr,epsfig,exscale,times,colordvi,
            latexsym,graphicx,subfigure}
\usepackage[english]{babel}
\usepackage{color}
\usepackage{anyfontsize}
\usepackage{float}
\usepackage{wrapfig}

\newcommand{\bm}[1]{\text{\boldmath $#1$\unboldmath}}  

\DeclareUnicodeCharacter{2061}{}

\newcommand{\vect}[1]{\mathbf{#1}}
\newcommand{\mat}[1]{\mathbf{#1}}

\newcommand{\nsd}{d}

\newcommand{\bn}{\bm{n}}
\newcommand{\bw}{\bm{w}}

\newcommand{\q}{\bm{q}}

\newcommand{\A}{\mat{A}}

\newcommand{\grad}{\bm{\nabla}}
\newcommand{\Div}{\bm{\nabla}\cdot}

\newcommand{\dO}{{\,dV}}
\newcommand{\dG}{{\,dS}}
\newcommand{\mean}[1]{\{#1\}}
\newcommand{\jump}[1]{\llbracket#1\rrbracket}

\newcommand{\hu}{\widehat{u}}

\newcommand{\traceh} {\ensuremath{\Lambda^h}}
\newcommand{\Vh} {\ensuremath{\mathcal{V}^h}}
\newcommand{\eltwo}{\ensuremath{\mathcal{L}^{2}}}
\newcommand{\Pki}{\ensuremath{\mathcal{P}_{k}(K_i)}}
\newcommand{\Pk}{\ensuremath{\mathcal{P}_{k}}}
\newcommand{\Aki}{\ensuremath{\mat{A}^{K_i}}}
\newcommand{\huIi}{\widetilde{u}^i}
\newcommand{\huI}{\widetilde{u}}

\newcommand{\bhuIi}{\widetilde{\vect{u}}^i}
\newcommand{\extBoundary}{\partial \Omega^{ext}}

\newcommand{\numel}{\ensuremath{{\texttt{n}_{\texttt{el}}}}}
\newcommand{\nOfElementFaces}{\ensuremath{n}}
\newcommand{\F}{\mat{F}}
\newcommand{\I}{\mathcal{I}}
\newcommand{\nface}{{\ensuremath{\texttt{n}_{\texttt{fc}}}}}

\begin{document}

\numberwithin{equation}{section}
\numberwithin{figure}{section}
\numberwithin{table}{section}

\newtheorem{theorem}{Theorem}[section]
\newtheorem{proposition}[theorem]{Proposition}
\newtheorem{lemma}[theorem]{Lemma}
\newtheorem{corollary}[theorem]{Corollary}

\newdefinition{remark}[theorem]{Remark}
\newproof{proof}{Proof}


\begin{frontmatter}

\title{eXtended Hybridizable Discontinous Galerkin (X-HDG) Method for Linear Convection-Diffusion Equations on Unfitted Domains}
\author[khas]{Haroon Ahmad}
\ead{haroon.ahmad@khas.edu.tr}
\author[khas]{Ceren G\"urkan\corref{cor1}}
\ead{ceren.gurkan@khas.edu.tr}
\cortext[cor1]{Corresponding author}
\address[khas]{Department of Civil Engineering, Kadir Has University, TR-34083 Cibali, Istanbul, Turkey}


\begin{abstract}
In this work, we propose a novel strategy for the numerical solution of linear convection diffusion equation (CDE) over unfitted domains. In the proposed numerical scheme, strategies from high order Hybridized Discontinuous Galerkin method and eXtended Finite Element method is combined with the level set definition of the boundaries. The proposed scheme and hence, is named as eXtended Hybridizable Discontinuous Galerkin (XHDG) method. In this regard, the Hybridizable Discontinuous Galerkin (HDG) method is eXtended to the unfitted domains; i.e, the computational mesh does not need to fit to the domain boundary; instead, the boundary is defined by a level set function and \textit{cuts} through the \textit{background mesh} arbitrarily. The original unknown structure of HDG and its hybrid nature ensuring the local conservation of fluxes is kept, while developing a modified bilinear form for the elements cut by the boundary. At every cut element, an auxiliary nodal trace variable on the boundary is introduced, which is eliminated afterwards while imposing the boundary conditions. Both stationary and time dependent CDEs are studied over a range of flow regimes from diffusion to convection dominated; using high order $(p\leq4)$ XHDG through benchmark numerical examples over arbitrary unfitted domains. Results proved that XHDG inherits optimal $(p+1)$ and super $(p+2)$ convergence properties of HDG while removing the fitting mesh restriction.    

\end{abstract}


\end{frontmatter}

\section{Introduction}
Linear convection diffusion (CD) equation can be of interest for all physical phenomena that involves a passive scalar being transported with a flowing medium. This can be interpreted as pollutant transport in a river or plumes carried away after a volcanic eruption. Numerical solution of convection diffusion equation and hence is of interest for scientists and engineers. Because of their discontinuous nature, convection dominated convection-diffusion problems; when numerically analyzed using classical continuous Galerkin (CG) methods, shows very poor stability properties and the numerical solution is often polluted by non-physical oscillations. For this reason CG methods, when used to solve convection dominated flow problems, often used together with a stabilization strategy. Among those stabilization strategies, the streamline upwind Petrov–Galerkin (SUPG) method \cite {Hughes-Brooks:79}, residual-free bubble method \cite{Brezzi-Russo:94} or adaptive mesh refinement techniques \cite{Eriksson-Estep-Hansbo-Johnson:95} are commonly used. For a through review on CDEs and their solution with stabilized continuous Galerkin methods the reader is referred to \cite{Morton:96, Roos-Stynes-Tobiska:96}. 

Discontinuous Galerkin (DG) methods on the other hand; are seen to be a more suitable alternative to solve CDEs especially because of their favorable stability properties. The earliest work in DG front is presented in \cite{Reed-Hill:73} for neutron transport and the recent developments can be found in \cite{Cockburn-Karniadakis-Shu:00}. Error estimates, convergence and accuracy of various DG methods for pure diffusion \cite {Becker-Hansbo-Larson:03, Karakashian-Pascal:03, Houston-Schotzau-Wihler:07}; for stationary CD \cite {Schotzau-Zhu:09,Ern-Stephansen-Vohralik:10} and for non stationary CD problems \cite{Ern-Proft:05, Houston-Suli:01, Georgoulis-Hall-Houston:08} are well understood in the community. Despite their advantageous stability properties, DG methods are often criticized because the number of degrees of freedom at their final system is more than its continuous counterpart for the same mesh and the same approximation degree.  

Hybridizable Discontinuous Galerkin (HDG) method \cite{Cockburn-CDG:08,Jay-CGL:09} is proposed as a possible remedy to that chokepoint. HDG inherits all advantageous properties of standard DG methods such as local conservation, built in stabilization, suitability for adaptivity and parallel computation while outperforming them with its hybrid and super convergent nature. Its hybrid nature reduces the number of degrees of freedom at the final system, similar to static condensation in CG methods \cite {GG-GMFH:13}. Moreover, HDG is based on a mixed formulation where both primal unknown and its derivative are approximated using same polynomial degree $p$ and this brings a $p+1$ convergence order in $\eltwo$ error norm not only for the primal unknown but also for its derivative. Consequently, a simple element by element post processing of the derivative leads to a super convergent approximation of order $p+2$ for the primal unknown \cite {Cockburn-CQS:12b}. HDG method is used to solve linear convection-diffusion equations successfully over standard fitted meshes in \cite {Nguyen-Peraire-Cockburn:09} and on general polyhedral fitted meshes in \cite{Qui-Sh:16}. 

The novelty of the strategy proposed here lays on without losing the favorable stability properties; optimal and super convergence and hybrid nature of the HDG method into an \textit{unfitted} setting; removing the mesh-fitting restriction and costs. 

\subsection{Earlier work} \label{sec:earlierwork}

The essential precondition for any numerical analysis is a high quality mesh resolving all the geometric features of the domain, ensuring an accurate enough numerical solution. Regardless of increasing computer power and ever growing research on accurate mesh generation; for complex stationary domains, mesh generation is still a challenge. For the problems involving complex domains or non-stationary interfaces; the majority of the computational cost is originated from accurate mesh creation or re-meshing. Independent of the numerical technique used, whenever a fitted mesh is on the table, scientists and engineers had to cope with these meshing/remeshing costs. To soothe this problem, so called, \textit{unfitted} discretization techniques are proposed. In an \textit{unfitted} setting, the mesh does not need to fit to the domain boundary, instead, a simple, cheap and stationary background mesh is created and the boundary or interface of interest is defined -generally- by a level set function \cite{Sethian} so that it can \textit{cut} through this background mesh. With unfitted discretizations, complex domains or non stationary interfaces can be studied avoiding the accurate fitted mesh creation or re-meshing costs. In the context of CG methods, several alternative unfitted/cut discretization techniques are proposed. Under the name of CutFEM \cite{Burman-Claus-Hansbo:15}, elliptic interface problems \cite{Burman-Zunino:12}, fluid structure interaction type coupled problems \cite{Massing-Larson-Logg-Rognes:15}, Stokes equations \cite{Burman-Claus-Massing:15} and Oseen equations \cite{Massing-Schott-Wall:18} are treated successfully. Unfitted discretizations are also used for the numerical solution of differential equations defined on surfaces \cite{Burman-Hansbo-Larson:15} or on surface-bulk coupled problems \cite{Hansbo-Larson-Zahedi:16}. 

Unfitted discretizations are as well adapted to DG framework to treat interface or evolving boundary problems, among many others
in \cite{Bastian-Engwer:09} scalar elliptic equation, in \cite{Bastian-Engwer-Fahlke-Ippisch:11} advection-diffusion equation, in \cite{Heimann-Engwer-Ippisch-Bastian:13}  incompressible Navier Stokes equations and in \cite{Gurkan-Sticko-Massing:20} advection-reaction problem with a high order CutDG method is analyzed. For the solution of surface and coupled surface-bulk PDEs using unfitted DG methods, the reader is referred to \cite{Massing:17} and the references therein. 

Building upon the HDG method and extending it to an unfitted setting; the XHDG method was initially proposed to treat void boundaries in \cite{Gurkan-Lardies-Kronbichler-Fernandez:16} and then implemented for the solution of elliptic bimaterial interface \cite{Gurkan-Kronbichler-Fernandez:17} and Stokes interface problems \cite{Gurkan-Kronbichler-Fernandez:18}, respectively.  

Animated by the efficiency of HDG on one hand and the effectiveness of unfitted discretizations on the other, the goal of novel XHDG method is to outreach the already existing analysis done over CDEs over an unfitted DG setting. There is a limited number of research done on the solution of CDEs using unfitted DG methods and the work presented here aims to outperform previous analysis from various aspects. The unfitted DG analysis done on CDEs in \cite{Bastian-Engwer-Fahlke-Ippisch:11} differs from the work presented here, since authors showed convergence only up to $p \leq 2$ while using second order approximation functions $p=2$. This means, the scheme employed was only $p$ order accurate when using $p$ order approximation, moreover, it did not have hybridization. Being different from steady unfitted HDG analysis for advection-reaction problems in \cite{Cockburn-Solano:14} and for steady Oseen equations in \cite{Solano-Vargas:22}; in this work with XHDG, no extrapolations or mesh-distance-to-boundary restrictions are at stake; leading to a more robust method, while studying both steady and time dependent problems. Furthermore, in both of the analysis, the authors showed convergence only up to $p \leq 3$. This work presents the novel XHDG method for the solution of CDEs, surpassing the previous work done in this field with its high order, unfitted, hybrid, robust and superconvergent nature with no mesh-distance limitations.   

\subsection{Contributions and outline of this paper}
In this work, we extend the XHDG discretization for the solution of steady and time dependent linear convection diffusion equations. Together with the level set description of the boundary, novel XHDG discretization for linear CD equation is proposed. Extension to time dependent CDEs is as well presented, using Backward Euler scheme for time discretization. At the elements \textit{cut} by the boundary, numerical integrations are handled using high order modified quadratures. Both diffusion and convection dominated regimes are studied while implementing two alternative built in flux stabilizations, namely, centered and upwind flux. When approximation functions of order $p$ is used, optimal $p+1$ convergence rate is reached both for primal unknown and its derivative; moreover, after an easy post process; super $p+2$ convergence rate is reached for the primal unknown.

This paper is organized as follows. The details of the XHDG discretization is presented in section \ref{sec:X-HDGFormulation}. In subsection \ref{sec:XHDG-localProblemNotCut}, local XHDG discretization of the elements that are not cut by the boundary is shown. For the elements not cut by the boundary, the element-wise formulation is equivalent to the standard HDG formulation. In subsection \ref{sec:XHDG-localProblemCut-Neumann}, the novel local weak form for the elements cut by the boundary when Neumann boundary conditions are imposed at the cut boundary  is presented. The local problem is then assembled into the so called global problem in subsection \ref{sec:globalProblem} where the hybridization happens by means of creating the system matrix in terms of the unknown only defined at element faces. Four numerical examples to prove the accuracy and convergence of XHDG is presented in section \ref{sec:numericalExamples}. For the first two examples the cut boundary is in simple circular shape; diffusion and convection dominated regimes examined respectively. For the third example a more complex "peanut" shaped cut boundary is studied. Lastly, the performance of XHDG for time dependent CDE is proved through the transport of a decaying Gaussian pulse. All numerical examples proved that XHDG reaches optimal and super convergence without computational mesh fitting to the domain boundary. This work is concluded with section \ref{sec:conclusions} stating the conclusions, final remarks and future work.  
%
\section{X-HDG formulation for linear convection-diffusion equations} \label{sec:X-HDGFormulation}
\begin{figure}[h]
\centering
\subfigure[]{
\includegraphics[width=0.42\textwidth]{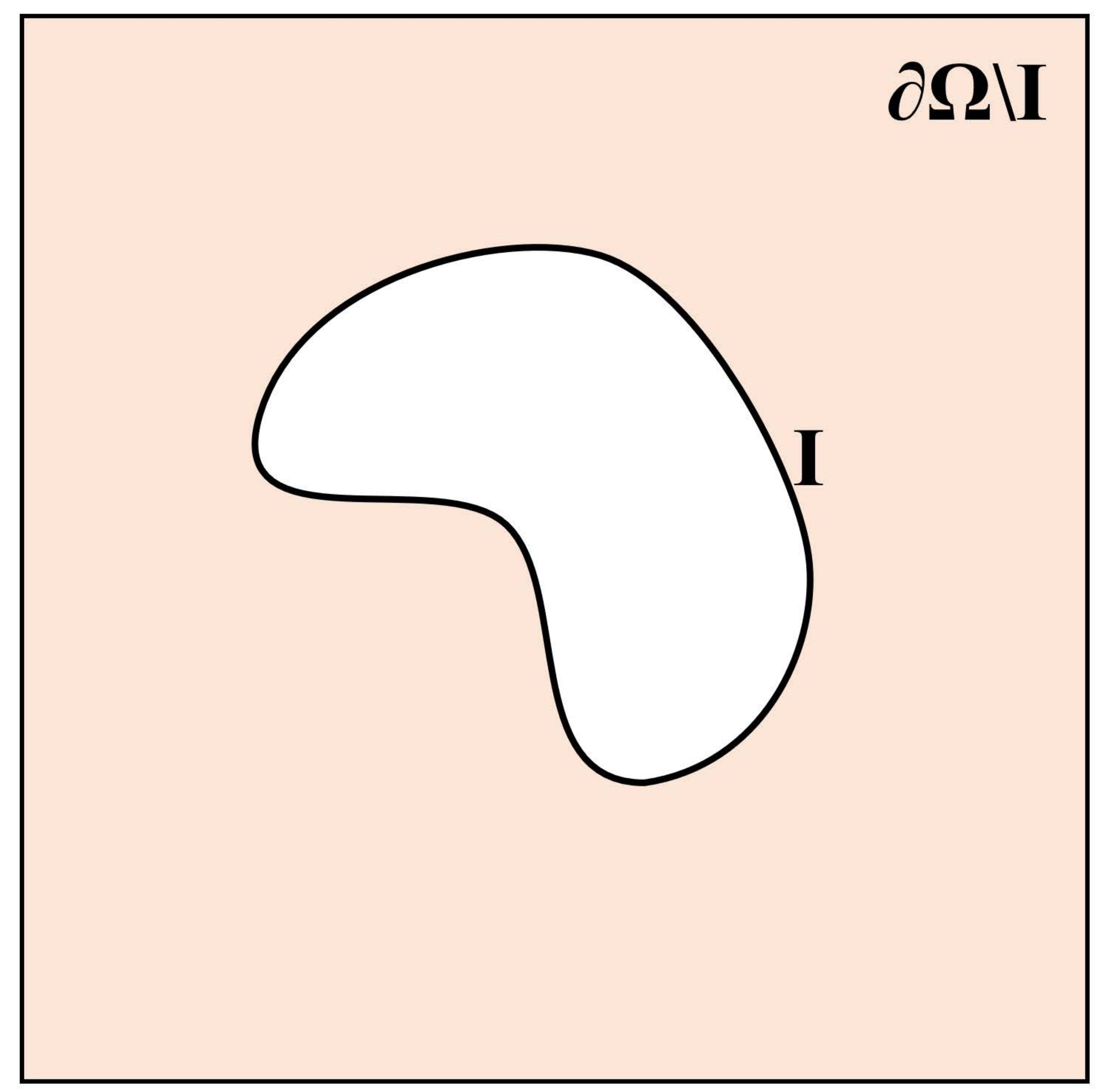}\label{fig:subfig1r}
}
\subfigure[]{
\includegraphics[width=0.42\textwidth]{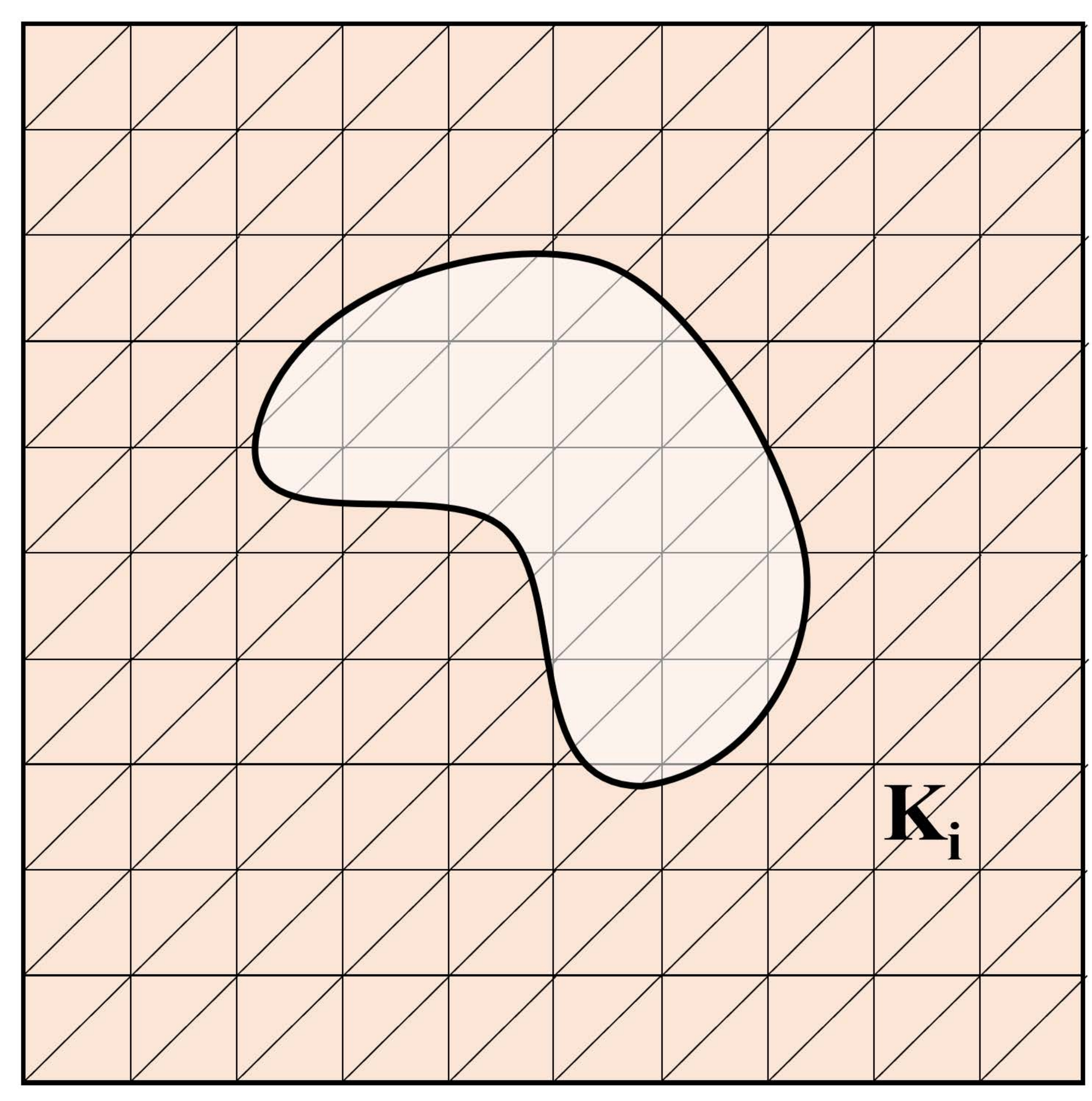}\label{fig:subfig2r}
}
\caption{The flow domain with an internal void: \subref{fig:subfig1r} schematic of the circular void inside a square domain $\Omega$ and \subref{fig:subfig2r} the background mesh fitted inside the domain. The void boundary $\I$ cutting the background mesh. Here, $\Omega$ is the overall domain, $\extBoundary$ is the external boundary, $\I$ is the internal void boundary and $K_{i}$ is a discrete triangular element.}
\label{fig:DomainVithCircularVoidExample}
\end{figure}

Consider a flow domain $\Omega$, bounded with an external boundary $\extBoundary:=\partial\Omega\backslash \I$ together with an arbitrary inner void boundary denoted by $\I$ as shown in Figure \ref{fig:DomainVithCircularVoidExample}. 
%
%
The steady, linear convection-diffusion model problem for this flow domain is defined as: 
%
\begin{equation}
\grad \cdot (\bm{c} u)-\grad \cdot (\nu \grad u) = f \quad \text{in } \Omega\label{eq:cd_eqn}
\end{equation}
%
where $u$ is the flow velocity, $c$ is a smooth velocity vector field, $\nu$ is kinematic viscosity and $f$ is the source term. The problem then can be closed with the following boundary conditions:
%
\begin{equation}
\begin{array}{rl}
u = u_D   \quad &\text{on } \extBoundary\\
(\bm{c}u+\bm{q})\cdot \bm{n} = g_N \quad &\text{on }  \I 
\end{array}\label{eq:cd_bc}
\end{equation}
where, $u_D$ is prescribed Dirichlet condition on the exterior boundary $\extBoundary$, and $g_N$ is a prescribed Neumann boundary condition on the inner boundary $\I$. 
The domain $\Omega$ is covered with a background finite element mesh with $\numel$ disjoint elements $K_i$, such that
%
\begin{equation*}
 \overline{\Omega}\subset \bigcup_{i=1}^{\numel}\overline{K}_i, \quad
 K_i\cap K_j=\emptyset ; \quad \text{for $i\neq j$,}
\quad  \extBoundary \subset \partial \left[ \bigcup_{i=1}^{\numel}\overline{K}_i \right]
\end{equation*}
%
Note that the background mesh fits the exterior boundary $\extBoundary$, but some elements may be cut by the interior boundary $\I$ as shown in Figure \ref{fig:DomainVithCircularVoidExample} right panel. The union of all $\nface$ faces $\Gamma_i$ (sides for 2D) intersecting the domain $\Omega$ is denoted as
%
\begin{equation*}
 \Gamma:=\bigcup_{i=1}^{\numel}\left[\partial K_i  \cap \overline\Omega\right]= \bigcup_{f=1}^{\nface}\left[\Gamma_f  \cap \overline\Omega\right].
\end{equation*}
%
The convection diffusion equation given in \eqref{eq:cd_eqn} in its strong form together with the boundary conditions is re-stated for every element $K_i$ that is not cut and cut by the inner boundary $\I$, following a mixed formulation in \eqref{eq:CDForHDG_standardElement} and \eqref{eq:CDForHDG_cutElement}, respectively. This element by element mixed formulation is so called as the \textit{local problem} in strong form.  
%
\begin{subequations}\label{eq:CDForHDG}
\begin{align}
\left.
\begin{array}{rl}
\bm{\nabla \cdot} (\bm{c}u + \bm{q}) = f \quad &\text{in } K_i \\
\bm{q} + \nu \bm{\nabla} u = 0 \quad &\text{in } K_i \\
u = \hu   \quad &\text{on } \partial K_i
\end{array}
\right\} & \text{ if } K_i \subset \Omega
\label{eq:CDForHDG_standardElement}
\\[0.5em]
\left.
\begin{array}{rl}
\bm{\nabla \cdot} (\bm{c}u + \bm{q}) = f  \quad &\text{in } \Omega_i  \\
 \bm{q} + \nu \bm{\nabla} u = 0 \quad &\text{in } \Omega_i\\
(\bm{c}u+\bm{q})\cdot \bm{n} = g_N  \quad &\text{on }  \I_i \\
u = \hu   \quad &\text{on } \partial \Omega_i \backslash\I_i
\end{array}
\right\} & \text{ if } \I \cap K_i \neq \emptyset
\label{eq:CDForHDG_cutElement}
\end{align}
\end{subequations}

for $i=1,\dotsc,\numel$, where, for cut elements,
\begin{equation}\label{eq:cutElement}
\Omega_i:=\Omega \cap K_i, \quad \I_i:=\I \cap K_i,
\end{equation}

\begin{wrapfigure}{l}{0.5\textwidth}
\begin{center}
\includegraphics[width=0.5\textwidth]{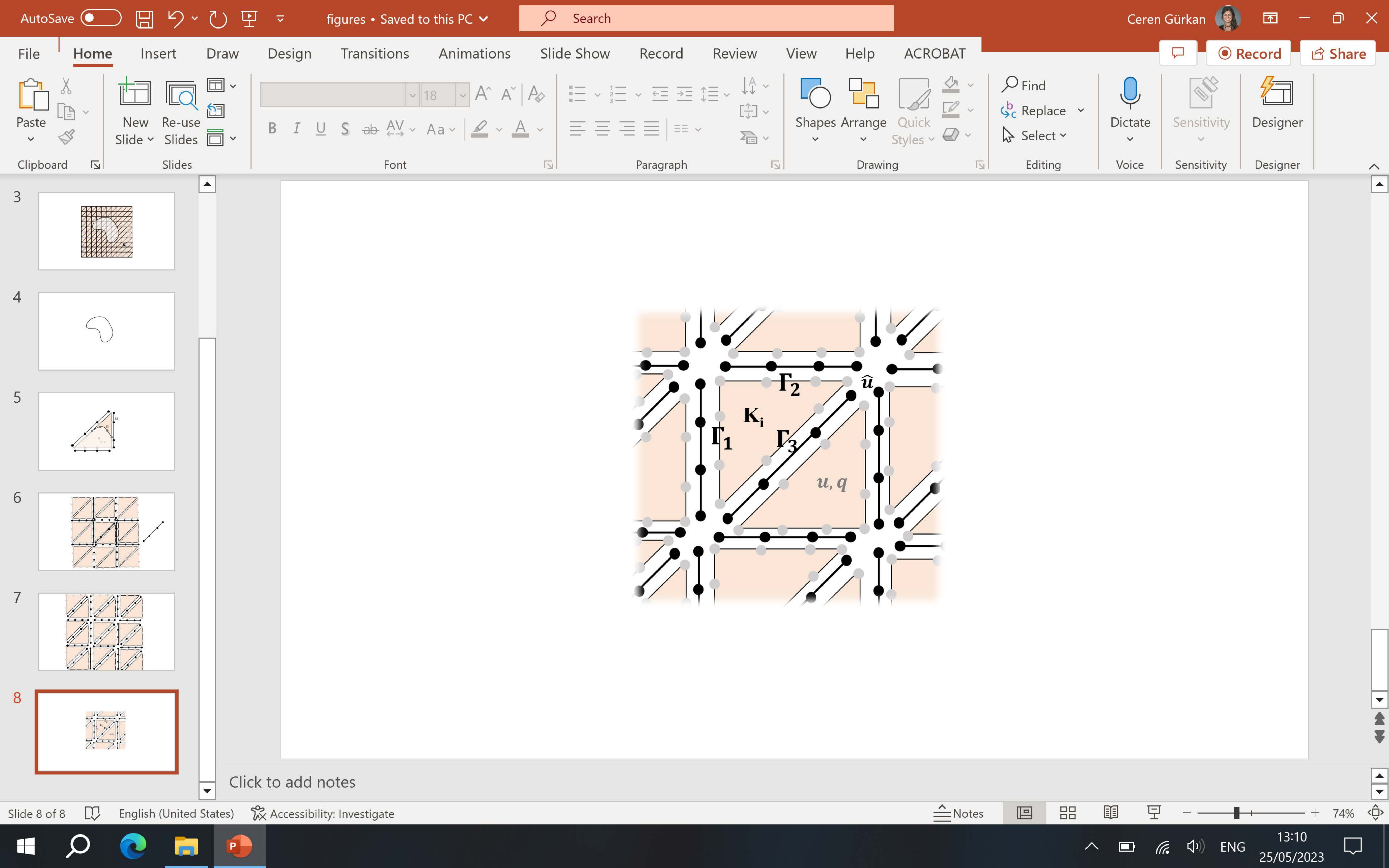}
\end{center}
\caption{Third degree XHDG discretization for a triangular standard element $K_i$ with three element faces $\Gamma_1$, $\Gamma_2$ and $\Gamma_3$.}
\label{fig:GlobalProb}
\end{wrapfigure}

%
In the local problem, two new variables, namely $\q$ and $\hu$ are introduced. The auxiliary variable, $\q$, corresponds to the flux of $u$, i.e., $\q:=\nu \bm \nabla u $. Introducing this flux variable and re-stating the \eqref{eq:CDForHDG} in a mixed form, leads to a numerical method that is $p+1$ convergent not only for the primal unknown $u$ but also for its derivative $\q$. The other new variable in the local problem is the $\hu$, which corresponds to the trace of $u$ at the element faces, see Figure \ref{fig:GlobalProb}. 

This trace variable is essential for the HDG and inherited in XHDG formulations since, the overall system matrix will include only this trace variable $\hu$ as an unknown. This leads to a significantly small system matrix and hence the name \textit{hybrid} is originated. It should be noted that, given the trace variable $\hu$, which is single valued at the element faces i.e., the mesh skeleton $\Gamma$, the local problems \eqref{eq:CDForHDG_standardElement} and \eqref{eq:CDForHDG_cutElement} can be solved at each element to determine the primal unknown $u$ and the flux $\q$.
Thus, the problem now reduces to determine the trace variable $\hu$ only, with the so-called \emph{conservativity conditions} (also known as \emph{global equations}), that is, the continuity of the flux across element faces, stated in \eqref{eq:conservativityEq}
\begin{equation} \label{eq:conservativityEq}
\jump{(\bm{c}u+\bm{q})\cdot \bm{n}} = 0 \quad \text{on } \Gamma \backslash \extBoundary,
\end{equation}
together with the exterior boundary condition
\begin{equation} \label{eq:DirichletCondition}
\hu =\mathbb{P}_2 (u_D) \quad \text{on } \extBoundary,
\end{equation}
the problem can be closed. 
The \emph{jump} $\jump{\cdot}$ operator is defined at a face $\Gamma_f$ as
\begin{equation*}
 \jump{\circledcirc} = \circledcirc_{L(f)} + \circledcirc_{R(f)} \quad \text{on } \Gamma_f,
\end{equation*}
%
where $R(f)$ and $L(f)$ are the left and the right elements sharing the face, that is, $\Gamma_f=\overline{K_{L(f)}}\cap \overline{K_{R(f)}}$, and the subindex $\circledcirc_{i}$ denotes the value of function $\circledcirc$ from element $K_i$. In particular,
\begin{equation*}
\jump{\q\cdot\bn}=\q_{L(f)}\cdot\bn_{L(f)}+\q_{R(f)}\cdot\bn_{R(f)}=(\q_{L(f)}-\q_{R(f)})\cdot\bn_{L(f)}
\end{equation*}
It is important noting that the continuity of the solution $u$ across $\Gamma$ is imposed by the Dirichlet boundary condition in the local problems $u=\hu$ and the fact that $\hu$ is single valued on $\Gamma$. The discretization of the conservativity condition \eqref{eq:conservativityEq} and the local problems \eqref{eq:CDForHDG}, with the boundary condition \eqref{eq:DirichletCondition}, leads to the complete X-HDG formulation.
The following discrete spaces for elemental variables, $u$ and $\q$, and for the trace variable, $\hu$, are considered
\begin{equation}\label{eq:discreteSpaces}
\begin{array}{l}
\Vh :=\bigl\{v \in\eltwo(\Omega) \; : v\vert_{K_i\cap \Omega} \in \Pk(K_i\cap \Omega)\; \text{ for }\; i=1,\dots,\numel \bigr\}
\\[1em]
W^h :=\bigl\{\bw \in\eltwo(\Omega) \; : \bw\vert_{K_i\cap \Omega} \in \Pk(K_i\cap \Omega)\; \text{ for }\; i=1,\dots,\numel \bigr\}
\\[1em]
\traceh :=\bigl\{\hat{v} \in\eltwo(\Gamma) \; : \hat{v}\vert_{\Gamma_f\cap \overline\Omega} \in \Pk(\Gamma_f\cap \overline\Omega) \; \text{ for }\; f=1,\dots,\nface \bigr\},
\end{array}
\end{equation}
where $\Pk$ denotes the space of polynomials of degree less or equal to $p$. In the coming subsections, the weak form of the local problem for standard and cut elements and the assembly of this elemental contributions using the global problem is presented.

\subsection{Local problem for standard elements} \label{sec:XHDG-localProblemNotCut}

In this subsection, the local XHDG formulation for an element $K_i$ not cut by the inner boundary $\I$; i.e., a standard element, is presented. For a standard element, XHDG formulation is identical to the HDG formulation. The local problem for a standard element in its strong form is given in \eqref{eq:CDForHDG_standardElement}. The weak form of the local problem for a standard element is stated in \eqref{eq:CDlocalproblemWeakStdElem}, that is: given $\hu\in \traceh$, find $u\in\Pki$, $\q\in[\Pki]^\nsd$ such that 
\begin{equation}
\begin{array}{lc}
\displaystyle -\int_{K_i} \bm{c}u \cdot \bm{\nabla}v \dO + \int_{\partial K_i} \tau (u -\hu) v \dG +
\displaystyle \int_{K_i}  v \Div \q \dO + \\
\hspace{10.5em} \displaystyle \int_{\partial K_i} (\bm{c}\cdot \bm{n}) \hu v \dG  = \displaystyle \int_{K_i}  vf \dO \quad \forall v\in \Pki 
\\
\displaystyle \int_{K_i} \dfrac{\q}{\nu}\cdot \bm{w} \dO  - \int_{K_i} u \Div \bm{w} \dO + \int_{\partial K_i} \hu (\bm{w}\cdot \bn) \dG = 0 \quad \forall \bm{w} \in [\Pki]^\nsd 
\label{eq:CDlocalproblemWeakStdElem}
\end{array}
\end{equation}
%
The first equation in \eqref{eq:CDlocalproblemWeakStdElem} can be derived from the first equation in \eqref{eq:CDForHDG_standardElement}, multiplying it with the test function $v$, applying integration by parts and replacing the flux by the numerical flux
%
\begin{equation}\label{eq:numericalFlux}
\widehat{\bm{c}u}+\bm{\hat{q}}:= \bm{c}\hat{u}+\bm{q}+\tau (u-\hat{u})\bm{\hat{n}},
\end{equation}
%
and undoing the integration by parts. The second equation, likewise, is obtained from the weak form of the second equation in \eqref{eq:CDForHDG_standardElement}, and applying integration by parts. In the weak form, whenever the primal variable $u$ is integrated over element face $\partial K_i$ it is replaced with the trace variable $u=\hu$.

The integral weak form of the local problem presented in \eqref{eq:CDlocalproblemWeakStdElem} can be re-written in matrix form as follows: 
\begin{equation} \label{eq:CDelementalSystemStdElem}
\left\{
\begin{array}{rl}
\left[\Aki_{uu} + \A^{\partial K_i}_{uu}\right] \vect{u}^i + \Aki_{uq} \vect{q}^i + \A^{\partial K_i}_{u\hu}  \vect{\Lambda}^{i}&= \vect{f}_{u}^{K_i}\\[0.5em]
\Aki_{qu} \vect{u}^i + \Aki_{qq} \vect{q}^i + \A^{\partial K_i}_{q\hu}  \vect{\Lambda}^{i}&= \vect{0}\\
\end{array}
\right.
\end{equation}
%
where $\vect{u}^i$ and $\vect{q}^i$ are the vectors of nodal values of $u$ and $\q$ in element $K_i$, and $ \Lambda^{i}$ is the vector of nodal values of $\hu$ on the $n$ faces of the element ($n=3$ for triangles and $n=4$ for tetrahedra). That is,
%
\begin{equation}
 \vect{\Lambda}^{i} := \left[
 \begin{array}{c}
 \widehat{\vect{u}}^{\F_{i1}}\\
 \vdots \\
 \widehat{\vect{u}}^{\F_{i\nOfElementFaces}}
 \end{array}
 \right],
\end{equation}
%
where $\widehat{\vect{u}}^f$ denotes the nodal values of $\hu$ on face $\Gamma_f$, and $\F_{ij}$ is the number of the $j$-th face of element $K_i$. 
Note that the subindeces in the $\A$ matrices refer to the space for the approximating function and the test function and superindices refer to the integral domain.

System \eqref{eq:CDelementalSystemStdElem} can be solved for $\vect{u}^i$ and $\vect{q}^i$ in each element, obtaining the so-called \emph{local solver} in the element $K_i$ such as:

\begin{equation}\label{eq:localSolver}
\vect{u}^i = \mat{U}^{K_i} \vect{\Lambda}^i+\vect{f}^{K_i}_U,\quad
\vect{q}^i = \mat{Q}^{K_i} \vect{\Lambda}^i+\vect{f}^{K_i}_Q,
\end{equation}
with
\begin{equation}\label{eq:localSolverEqStandardHDG}
\left[ \begin{array}{c}
\mat{U}^{K_i}\\[0.5em] \mat{Q}^{K_i}
\end{array}\right]
= -
 \mat{\mathbb{A}}^{-1}
\left[ \begin{array}{c}
 \Aki_{u\hu}\\[0.5em]  \Aki_{q\hu}
\end{array}\right],
  \;\;
\left[ \begin{array}{c}
\mat{f}^{K_i}_U\\[0.5em] \mat{f}^{K_i}_Q
\end{array}\right]
=
 \mat{\mathbb{A}}^{-1}
\left[ \begin{array}{c}
 \vect{f}_{u}^{K_i}\\[0.5em] \vect{0}
\end{array}\right]
\end{equation}
where
\begin{equation*}
 \mat{\mathbb{A}}= \left[ \begin{array}{cc}
\left[\Aki_{uu}+\A^{\partial K_i}_{uu}\right]  & \Aki_{uq} \\[0.5em] \Aki_{qu}  & \Aki_{qq}
\end{array}\right]
\end{equation*}
%
The local solver leads to a system where the elemental unknowns, namely; $\vect{u}^i$ and $\vect{q}^i$, can be explicitly expressed only in terms of the trace variable on element faces, $\vect{\Lambda}^i$.
%
\subsection{Local problem for a cut element-Neumann boundary conditions} \label{sec:XHDG-localProblemCut-Neumann}

In this subsection the XHDG local problem at an element $\Omega_i$ cut by the interior boundary $\I$ is detailed. Neumann boundary conditions are imposed at the inner boundary $\I$, and at the outer boundary $\partial\Omega^{ext}$, Dirichlet boundary conditions are considered. The local problem then can be defined as: given $\hu\in \traceh$, find $u\in\Pki$, $\q\in[\Pki]^\nsd$ such that  
%
\begin{equation}\label{eq:CDXHDGlocalproblemCutNeu}
\begin{array}{cl}
\displaystyle -\int_{\Omega_{i}} \bm{c}u\cdot \bm{\nabla}v \dO + \int_{\Omega_i}  v \Div \q \dO + \int_{\partial \Omega_{i}\backslash \I_{i}} \tau (u-\hu) v\dG + \int_{\partial \Omega_{i}\backslash \I_{i}}(\bm{c}\cdot \bm{n}) \hu v\dG \\
 + \displaystyle \int_{I_{i}} \tau (u - \huI) v \dG + \int_{I_{i}} (\bm{c}\cdot \bm{n})\huI v \dG  = \displaystyle \int_{\Omega_i}  vf \dO \qquad \forall v\in \Pki
\\[1em]
\displaystyle \int_{\Omega_i} \frac{\bm{q}}{\nu}\cdot \bm{w} \dO-\int_{\Omega_i} u \Div \bm{w} \dO + \int_{\partial \Omega_i\backslash \I_{i}} \hu (\bm{w}\cdot \bn) \dG \\
\qquad +\displaystyle \int_{\I_{i}} \huI (\bm{w}\cdot\bn) \dG 
= 0 \qquad \forall \bm{w} \in [\Pki]^\nsd 
\end{array}
\end{equation}
\newpage
%
\begin{wrapfigure}{l}{0.4\textwidth}
\begin{center}
\includegraphics[width=0.35\textwidth]{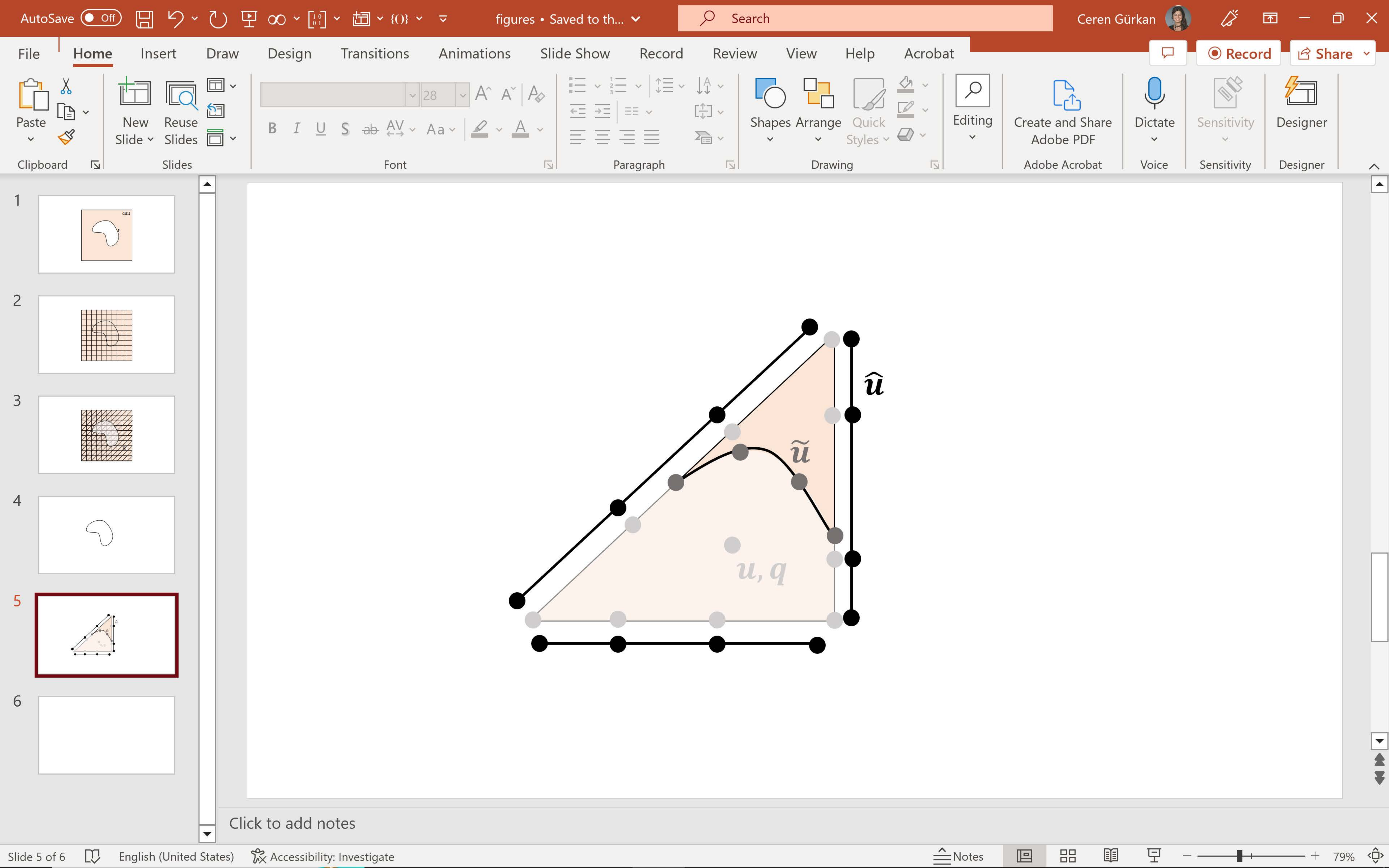}
\end{center}
\caption{Cut element with elemental variables $u$ and $\q$, trace variable $\hu$ and the variable $\huIi$ defined on inner boundary.}
\label{fig:CutElem}
\end{wrapfigure}
%
for all $v\in \Pk(\Omega_i)$ and $\bw\in [\Pk(\Omega_i)]^\nsd$, where $\Omega_i=\Omega \cap K_i$. In \eqref{eq:CDXHDGlocalproblemCutNeu} yet another unknown variable $\huI$  is defined on the interior boundary $\I_i=\I \cap K_i$. Compared to the weak form for standard elements \eqref{eq:CDlocalproblemWeakStdElem}, the XHDG weak form for a cut element has additional terms corresponding to integrals along the interface $\I_i$, involving the new trace variable $\huI$. This is because in a standard triangular element the surface integrals defined on $\partial K_i$ are the element faces only, however; in a cut element, not only the element faces but also an additional inner boundary forms the cut element, see Figure \ref{fig:CutElem}. 
The integral weak form of the local problem for a cut element given in \eqref{eq:CDXHDGlocalproblemCutNeu} can be restated in matrix form as follows: 

\begin{equation} \label{eq:CDMatrixSystemCutNeu}
\begin{array}{rl}
\left[\A^{\Omega_i}_{uu} + \A^{\partial \Omega_i \backslash \I_i}_{uu} + \A^{\I_i}_{uu}\right]\vect{u}^i &+ \A^{\Omega_i}_{uq} \vect{q}^i + \A^{\partial\Omega_i\backslash\I_i}_{u\hu} \vect{\Lambda}^{i} + \A^{\I_i}_{u\huI} \vect{\huI}^{i} = \vect{f}_u^{\Omega_i}\\[0.5em]
\A^{\Omega_i}_{qu} \vect{u}^i &+ \A^{\Omega_i}_{qq} \vect{q}^i + \A^{\partial \Omega_i \backslash \I_i}_{q\hu}  \vect{\Lambda}^{i}+
\A^{\I_i}_{q\huI}  \vect{\huI}^{i} = \vect{0}\\
\end{array}
\end{equation}
%
similar to \eqref{eq:CDelementalSystemStdElem}, but with three new matrices corresponding to integrals on the inner boundary $\I_i$, and the nodal values for the new trace variable, $\bhuIi$. To be able to obtain the \textit{local solver} similar to \eqref{eq:localSolver} we need to eliminate the new unknown variable $\bhuIi$ from \eqref{eq:CDMatrixSystemCutNeu}. By imposing the Neumann boundary condition at the inner boundary $\I_i$, $\bhuIi$ can be eliminated from \eqref{eq:CDMatrixSystemCutNeu}.     
The Neumann boundary condition at the inner boundary is stated as:

\begin{equation*}
(\widehat{\bm{c}u}+\widehat{\bm{q}})\cdot \bm{n} = g_N  \quad \text{on }  \I_i 
\end{equation*}
Replacing the numerical flux with the expression given in \eqref{eq:numericalFlux}, the weak form of the Neumann condition can be written as: given $u\in\Pk(\Omega_i)$, $\q\in[\Pk(\Omega_i)]^\nsd$, find $\huIi\in\Pk(\I_i)$ such that
\begin{equation}\label{eq:conservativity_Ii}
\int_{\I_i} (\bm{c}\widetilde u \cdot \bn) \widetilde v   + \int_{\I_i} (\q \cdot \bn) \widetilde v \dG + \tau \int_{\I_i} \widetilde v\left( u-\huI\right) \dG =  \int_{\I_i} \widetilde v g_N \dG
    \quad \forall \widetilde v\in  \Pk(\I_i)
\end{equation}
%
Equation \eqref{eq:conservativity_Ii} can be written in matrix form as follows
\begin{equation}
    \A^{\I_i}_{\huIi u}\vect{u}^i + \A^{\I_i}_{\huIi q}\vect{q}^i + \A^{\I_i}_{\huIi \huIi} \bhuIi =  \vect{g_N},
    \label{eq:NeumannCondMatrixForm}
\end{equation}
Using \eqref{eq:NeumannCondMatrixForm} the new variable $\huIi$ defined on inner boundary $\I$ can be written in terms of $u$ and $\bm{q}$ as follows

\begin{equation}\label{eq:conservativity_Ii_discrete}
    \bhuIi = \mat{T}^i_u\vect{u}^i + \mat{T}^i_q \vect{q}^i - \vect{t}^i.
\end{equation}
where,
\[
\mat{T}^i_u= - [\A^{\I_i}_{\huIi \huIi}]^{-1} \A^{\I_i}_{\huIi u},\quad \mat{T}^i_q = -[\A^{\I_i}_{\huIi \huIi}]^{-1} \A^{\I_i}_{\huIi q},
\quad \vect{t}^i = -[\A^{\I_i}_{\huIi \huIi}]^{-1}\vect{g_N}.
\] 
If in the local problem \eqref{eq:CDMatrixSystemCutNeu} $\bhuIi$ is replaced with the expression \eqref{eq:conservativity_Ii_discrete} the updated local problem in matrix form can be written as follows
\begin{equation*} \label{eq:elementalSystemCut2}
\begin{array}{ll}
\left[\A^{\Omega_i}_{uu} +\A^{\partial \Omega_i\backslash \I_i}_{uu} + \A^{\I_i}_{uu} +\A^{\I_i}_{u\huIi} \mat{T}^i_u\right]\vect{u}^i 
+ \left[\A^{\Omega_i}_{uq} +\A^{\I_i}_{u\huIi} \mat{T}^i_q\right] \vect{q}^i + \A^{\partial \Omega_i\backslash \I_i}_{u\hu}  \vect{\Lambda}^{i} \\[0.5em]
\hspace{8.51cm}= \vect{f}_u^{\Omega_i}+\A^{\I_i}_{u\huIi}\vect{t}^i 
\\[0.5em]
\left[\A^{\Omega_i}_{qu} + \A^{\I_i}_{q\huIi}  \mat{T}^i_u\right]\vect{u}^i
+ \left[\A^{\Omega_i}_{qq} + \A^{\I_i}_{q\huIi} \mat{T}^i_q\right]\vect{q}^i + \A^{\partial \Omega_i\backslash \I_i}_{q\hu}  \vect{\Lambda}^{i} = \A^{\I_i}_{q\huIi}\vect{t}^i\\
\end{array}
\end{equation*}
After eliminating the new variable $\huIi$ from local problem using Neumann boundary condition; identical to \eqref{eq:localSolver}, \textit{local solver} for a cut element can be written as:
\begin{equation}\label{eq:localSolverCut2}
\vect{u}^i = \mat{U}^{K_i} \vect{\Lambda}^i+\vect{f}^{K_i}_U,\quad
\vect{q}^i = \mat{Q}^{K_i} \vect{\Lambda}^i+\vect{f}^{K_i}_Q,
\end{equation}
with
\begin{equation} \label{eq:localSolverCut}
\left[ \begin{array}{c}
\mat{U}^{K_i}\\[0.5em] \mat{Q}^{K_i}
\end{array}\right]
= -
\mat{\mathbb{A}}^{-1}
\left[ \begin{array}{c}
 \A^{\partial \Omega_i \backslash \I_i}_{u\hu}\\[0.5em]  \A^{\partial\Omega_i\backslash \I_i}_{q\hu}
\end{array}\right],
  \;\;
\left[ \begin{array}{c}
\mat{f}^{K_i}_U\\[0.5em] \mat{f}^{K_i}_Q
\end{array}\right]
=
\mat{\mathbb{A}}^{-1}
\left[ \begin{array}{c}
 \vect{f}_{u}^{\Omega_i}+\A^{\I_i}_{u\huIi}\vect{t}^i
 \\[0.5em]
 \A^{\I_i}_{q\huIi}\vect{t}^i
\end{array}\right]
\end{equation}
and
\begin{equation*}
 \mat{\mathbb{A}}= \left[ \begin{array}{cc}
\left[\A^{\Omega_i}_{uu} +\A^{\partial \Omega_i\backslash \I_i}_{uu}+ \A^{\I_i}_{uu} +\A^{\I_i}_{u\huIi} \mat{T}^i_u\right]  & \left[\A^{\Omega_i}_{uq} +\A^{\I_i}_{u\huIi} \mat{T}^i_q\right] \\[0.5em]
\left[\A^{\Omega_i}_{qu} + \A^{\I_i}_{q\huIi}  \mat{T}^i_u\right]  & \left[\A^{\Omega_i}_{qq} + \A^{\I_i}_{q\huIi} \mat{T}^i_q\right]
\end{array}\right]
\end{equation*}
It is important to understand that the structure of the local solver is similar to the one for non-cut elements \eqref{eq:localSolver}. In the local solver unknown variables $u$ and $\bm{q}$ are once more expressed only in terms of the trace unknown $\hu$, which is the essence of hybrid discretization.  


\begin{remark}
The main difference between unfitted XHDG discretization and fitted HDG discretization lays in the local problem for cut elements. To calculate the partial integrals over cut elements and cut faces, a high order modified quadrature is introduced. The details of this implementation can be found in \cite{Gurkan-Lardies-Kronbichler-Fernandez:16}, section 4.   
\end{remark}


\begin{remark}
The implementation of Dirichlet boundary conditions on internal boundary is straightforward. Starting from the mixed strong form:
\begin{equation*}
\begin{array}{rl}
\bm{\nabla \cdot} (\bm{c}u + \bm{q}) = f  \quad &\text{in } \Omega_i  \\
 \bm{q} + \nu \bm{\nabla} u = 0 \quad &\text{in } \Omega_i\\
u = u_I \quad &\text{on }  \I_i \\
u = \hu   \quad &\text{on } \partial \Omega_i \backslash\I_i
\end{array}
\end{equation*}
And multiplying it with test functions, doing integration by parts, replacing the flux with the numerical flux $\widehat{\bm{c}u}+\bm{\hat{q}}:= \bm{c}(\hat{u}+u_I)+\bm{q}+\tau (u-\hat{u})\bm{\hat{n}}+\tau (u-u_I)\bm{\hat{n}}$ and undoing integration by parts (only for first equation); one can obtain the weak form  below for every cut element when Dirichlet boundary conditions are imposed at the cut boundary.  
\begin{equation*}
\begin{array}{cl}
\displaystyle -\int_{\Omega_{i}} \bm{c}u\cdot \bm{\nabla}v \dO + \int_{\Omega_i}  v \Div \q \dO + \int_{\partial \Omega_{i}\backslash \I_{i}} \tau (u-\hu) v\dG + \int_{\partial \Omega_{i}\backslash \I_{i}}(\bm{c}\cdot \bm{n}) \hu v\dG \\
+ \displaystyle \int_{I_{i}} \tau u v \dG = \displaystyle \int_{\Omega_i}  vf \dO - \displaystyle \int_{I_{i}} \tau u_I v \dG - \int_{I_{i}} (\bm{c}\cdot \bm{n})u_I v \dG  \qquad \forall v\in \Pki
\\[1em]
\displaystyle \int_{\Omega_i} \frac{\bm{q}}{\nu}\cdot \bm{w} \dO-\int_{\Omega_i} u \Div \bm{w} \dO + \int_{\partial \Omega_i\backslash \I_{i}} \hu (\bm{w}\cdot \bn) \dG \\
\qquad  
= -\displaystyle \int_{\I_{i}} u_I (\bm{w}\cdot\bn) \dG \qquad \forall \bm{w} \in [\Pki]^\nsd 
\end{array}
\end{equation*}
By close inspection, the reader can realize that the integral weak form given above for Dirichlet boundaries are quite similar to the one of Neumann boundaries; the only difference is that now, the integrals over cut boundary $\I$ with unknown $\huI$ do appear at the right hand side of the equations with known nodal Dirichlet values $u_I$.    
 
\end{remark}


\begin{remark}
The built in stabilization of XHDG is inherited from HDG. The built in stabilization is embedded in the numerical flux definition \eqref{eq:numericalFlux}. In \cite{Nguyen-Peraire-Cockburn:09} it is shown that the equations \eqref{eq:CDlocalproblemWeakStdElem} and \eqref{eq:CDXHDGlocalproblemCutNeu} together with the global problem have a unique solution if $\grad \cdot \bm{c} \geq 0$ and when stabilization parameter $\tau$ satisfies the condition $\tau > \frac{1}{2} \hspace{0.1cm} (\bm{c} \cdot \bm{n})$. 
The stabilization parameter $\tau$ here is chosen following the strategy presented in \cite{Nguyen-Peraire-Cockburn:09}. There are two possible choices, one leading to \textit{centered} and the other to \textit{upwind} flux definition. First, lets define $\tau=: \tau_d + \tau_c$ where $\tau_d$ represents diffusive and $\tau_c$ represents convective flux. When $\tau_d^+ = \tau_d^- = \eta_d$ and $\tau_c^+ = \tau_c^- = \eta_c$ with $\eta_d = \frac{\nu}{\textit{l}} $ and $\eta_c = |\bm{c} \cdot \bm{n}| $ the flux scheme is named as \textit{centered} where $\textit{l}$ denotes the diffusive length scale. In the centered flux definition we have a single-valued $\tau$, that is $\tau^+=\tau^-$. On the other hand, when we choose $\tau_d^{\pm}$ and $\tau_c^{\pm}$ as
\begin{equation*}
(\tau_d^{\pm}, \tau_c^{\pm}) = (\eta_d,\eta_c) \frac{|\bm{c} \cdot \bm{n}^+| \pm \bm{c} \cdot \bm{n}^+}{2|\bm{c} \cdot \bm{n}^+|}
\end{equation*}
the scheme leads to the \textit{upwind} flux definition. In the numerical examples, both centered and upwind flux definitions are used.  
\end{remark}

\subsection{Global problem} \label{sec:globalProblem}

Element by element contributions both from standard and cut elements, leads to the local solver \eqref{eq:localSolver} and \eqref{eq:localSolverCut2}. The local solver defines the unknowns in the element, $u$ and $\q$, in terms of the trace unknown at element boundary, $\hu$.
Therefore, the problem can be reduced to determine the trace unknowns $\{\widehat{\vect{u}}^f\}_{f=1}^{\nface}$ on the mesh skeleton $\Gamma$. To create the system matrix, element by element contributions are assembled using the so-called \emph{global problem}. Global problem corresponds to the discretization of the conservativity of flux on $\Gamma$, stated in \eqref{eq:conservativityEq}.

Replacing flux by the numerical flux \eqref{eq:numericalFlux}, the weak form of the conservativity condition can be stated as: find $\hu\in \traceh$ such that $\hu = u_D$ on $\extBoundary$ and
\begin{equation*}
    \int_\Gamma \widehat{v} \jump{\q \cdot \bn} \dG + 2\int_\Gamma \widehat{v}  \mean{\bm{c}.\bm{n}}\hu \dG + 2\int_\Gamma \widehat{v} (\mean{\tau u}-\mean{\tau}\hu) \dG = 0
    \quad \forall \widehat{v}\in  \traceh,
\end{equation*}
where $\mean{\cdot}$ is the mean operator on the faces defined as,
\begin{equation*}
    \mean{\circledcirc} = \frac{1}{2} \left( \circledcirc_{L(f)} + \circledcirc_{R(f)} \right) \quad \text{on }\Gamma_f.
\end{equation*}
Discrete form of the conservativity condition leads to the following matrix system
\begin{equation}\label{eq:CDGlobalProbMatrixForm}
    \A^{f,L}_{\hu u} \vect{u}^{L(f)} +\A^{f,L}_{\hu q} \vect{q}^{L(f)}
    \;+\; \A^{f,R}_{\hu u} \vect{u}^{R(f)} +\A^{f,R}_{\hu q} \vect{q}^{R(f)} \;+\; \A^f_{\hu\hu} \widehat{\vect{u}}^f= 0.
\end{equation}

Elemental contributions from both left $\vect{u}^{L(f)} \text{, } \vect{q}^{L(f)}$ and right $\vect{u}^{R(f)} \text{, } \vect{q}^{R(f)}$ side of the element face are assembled into \eqref{eq:CDGlobalProbMatrixForm} using the local solver. In \eqref{eq:CDGlobalProbMatrixForm} when $\vect{u}$ and $\vect{q}$ are replaced with local solver from standard \eqref{eq:localSolver} or from cut element \eqref{eq:localSolverCut2}, the final system then can be obtained only in terms of the trace variable $\vect{\Lambda}^i$.  

XHDG implementation, being similar to HDG implementation; involves a loop over elements. For each element, the matrices and vectors of the local solver \eqref{eq:localSolver} and \eqref{eq:localSolverCut2} are computed, and then those contributions are assembled into the system matrix using the global solver \eqref{eq:CDGlobalProbMatrixForm}. 
Once the system is assembled the Dirichlet boundary conditions on the outer boundary \eqref{eq:DirichletCondition} (when applicable) are imposed and then the system can be solved. Since in the global system the only unknown variable is trace variable $\vect\Lambda$ defined on mesh skeleton $\Gamma$; with the solution of this system, the nodal values of the trace variable is obtained. After the nodal trace values are obtained, with a simple put-back operation using the local solver, the elemental unknowns $\vect{u}^i$ and $\q^i$, can be computed.
It is important to note that in X-HDG implementation standard HDG structure is kept. The basic differences are the derivation of modified weak forms over cut elements, modified integral calculations at cut elements, the new unknowns defined at cut boundary $\huI$ and their elimination.
\begin{remark}
For the simplest XHDG implementation; all faces are assembled in the global problem even if they do not intersect the domain i.e. when they are inside the void part of the domain. The matrices assembled for the faces not intersecting the domain lead to null rows and columns in the global system. Those rows and columns are then removed from the system, reducing its size and obtaining a system with unique solution, together with imposing the Dirichlet boundary condition \eqref{eq:DirichletCondition}.
\end{remark}
\begin{remark} \label{rm:superconvergencePostprocess}
An element-by-element postprocessing can be performed to compute a superconvergent numerical solution from original numerical solution obtained after solving the global system. The superconvergence is one of the most distinguished properties of HDG, making it outperform other DG methods. This supercovergence is as well inherited by XHDG. Similar to standard HDG; the superconvergent solution can be computed in every element $K_i$ to find $u^*\in\mathcal{P}_{k+1}(\Omega_i)$ such that:
\begin{equation*}
\begin{array}{c}
\displaystyle\int_{\Omega_i} \nu \grad u^* \cdot \grad v \dO = -\int_{\Omega_i} \q \cdot \grad v \dO \quad \forall v \in \mathcal{P}_{k+1}(\Omega_i),\\
\displaystyle\int_{\Omega_i} u^* \dO = \int_{\Omega_i} u \dO,
\end{array}
\end{equation*}
with $\Omega_i = K_i$ for standard elements, and $\Omega_i = K_i\cap \Omega$ for cut elements. The solution of this element-by-element computation, $u^*$, converges with order $p+2$ in the $\mathcal{L}_2$ norm. See \cite{Cockburn-CDG:08,Cockburn-CQS:12b} for details and other possible computations of a superconvergent solution.
\end{remark}

\subsection{Extension to time dependent problems}
The extension of steady XHDG formulation to unsteady problems is straightforward. Here, the unsteady XHDG formulation is presented for a cut element, the unsteady weak form for a standard element can be derived following the same steps. The mixed strong form for a Neumann cut element can be defined as:
\begin{equation*}
\begin{array}{rl}
\dfrac{du}{dt}+\bm{\nabla \cdot} (\bm{c}u + \bm{q}) = f  \quad &\text{in } \Omega_i  \\
 \bm{q} + \nu \bm{\nabla} u = 0 \quad &\text{in } \Omega_i\\
(\bm{c}u+\bm{q})\cdot \bm{n} = g_N  \quad &\text{on }  \I_i \\
u = \hu   \quad &\text{on } \partial \Omega_i \backslash\I_i
\end{array}
\end{equation*}
when time dependent term is discretized using Backward Euler method:
\begin{equation*}
\displaystyle \int_{\Omega_i}\dfrac{du}{dt} v \dO = \displaystyle \int_{\Omega_i} \dfrac{1}{\Delta t} u^t v \dO - \displaystyle \int_{\Omega_i} \dfrac{1}{\Delta t} u^{t-1} v \dO \qquad \forall v\in \Pk(\Omega_i)
\end{equation*}
 The local unsteady weak form for a Neumann cut element can be stated as follows with two additional terms coming from time discretization: 
\begin{equation*}
\begin{array}{cl}
\displaystyle \int_{\Omega_i} \dfrac{1}{\Delta t} u^t v \dO \displaystyle -\int_{\Omega_{i}} \bm{c}u\cdot \bm{\nabla}v \dO + \int_{\Omega_i}  v \Div \q \dO + \int_{\partial \Omega_{i}\backslash \I_{i}} \tau (u-\hu) v\dG +\\ \displaystyle \int_{\partial \Omega_{i}\backslash \I_{i}}(\bm{c}\cdot \bm{n}) \hu v\dG 
 + \displaystyle \int_{I_{i}} \tau (u - \huI) v \dG + \int_{I_{i}} (\bm{c}\cdot \bm{n})\huI v \dG \\
\hspace{1.5cm} = \displaystyle \int_{\Omega_i} \dfrac{1}{\Delta t} u^{t-1} v \dO +\displaystyle \int_{\Omega_i}  vf \dO \qquad \forall v\in \Pki
\\[1em] \vspace{0.3cm}
\displaystyle \int_{\Omega_i} \frac{\bm{q}}{\nu}\cdot \bm{w} \dO-\int_{\Omega_i} u \Div \bm{w} \dO + \int_{\partial \Omega_i\backslash \I_{i}} \hu (\bm{w}\cdot \bn) \dG \\
+\displaystyle \int_{\I_{i}} \huI (\bm{w}\cdot\bn) \dG 
= 0 \qquad \forall \bm{w} \in [\Pki]^\nsd 
\end{array}
\end{equation*}
%
the local solver then can be obtained as explained in subsection \ref{sec:XHDG-localProblemCut-Neumann} and the global problem can be formed as stated in subsection \ref{sec:globalProblem}. 
\section{Numerical tests}\label{sec:numericalExamples}

The performance of X-HDG method is tested over four numerical examples where the boundary is defined by a level set function, cutting through the background mesh. In the first example, square background geometry with circular void boundary is considered on diffusion domainated flow regime. In the second example, the geometry and the boundary kept the same while changing the flow regime to be convection dominated. For the third example, an arbitrary cut boundary in peanut shape is solved on convection dominated regime; to demonstrate the capability of XHDG to handle arbitrary boundary cuts. Last but not least, transport of a decaying Gaussian pulse is studied. All tests prove that removing the boundary-fitting-mesh restriction and using a cut mesh did not effect the optimal and super convergence or accuracy of the XHDG when compared to its mesh-fitting counterpart, HDG.   

\subsection{Circular cut boundary on a square domain - Diffusion dominated}
\begin{figure}[h]
\centering
\subfigure[]{
\includegraphics[width=0.45\textwidth]{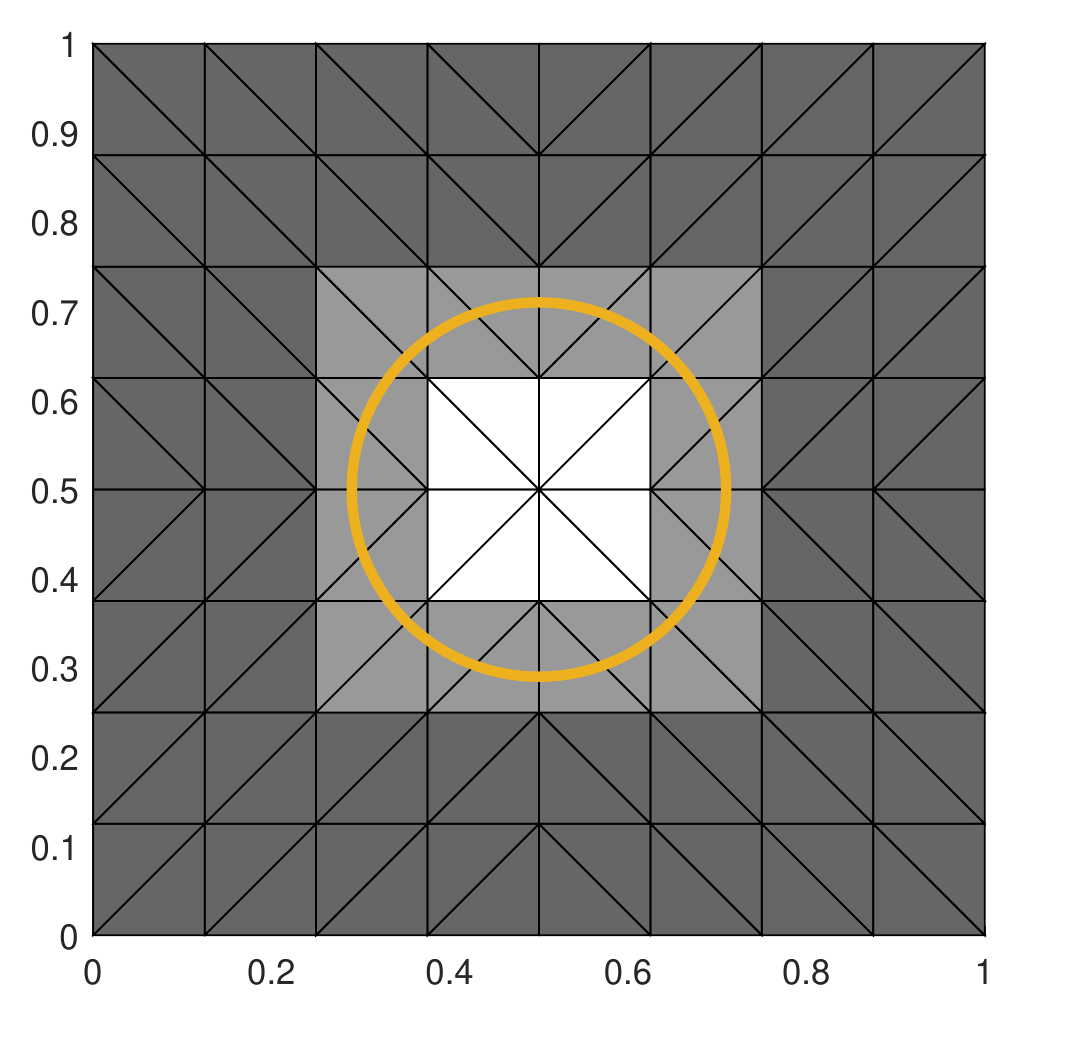}
}
\subfigure[]{
\includegraphics[width=0.44\textwidth]{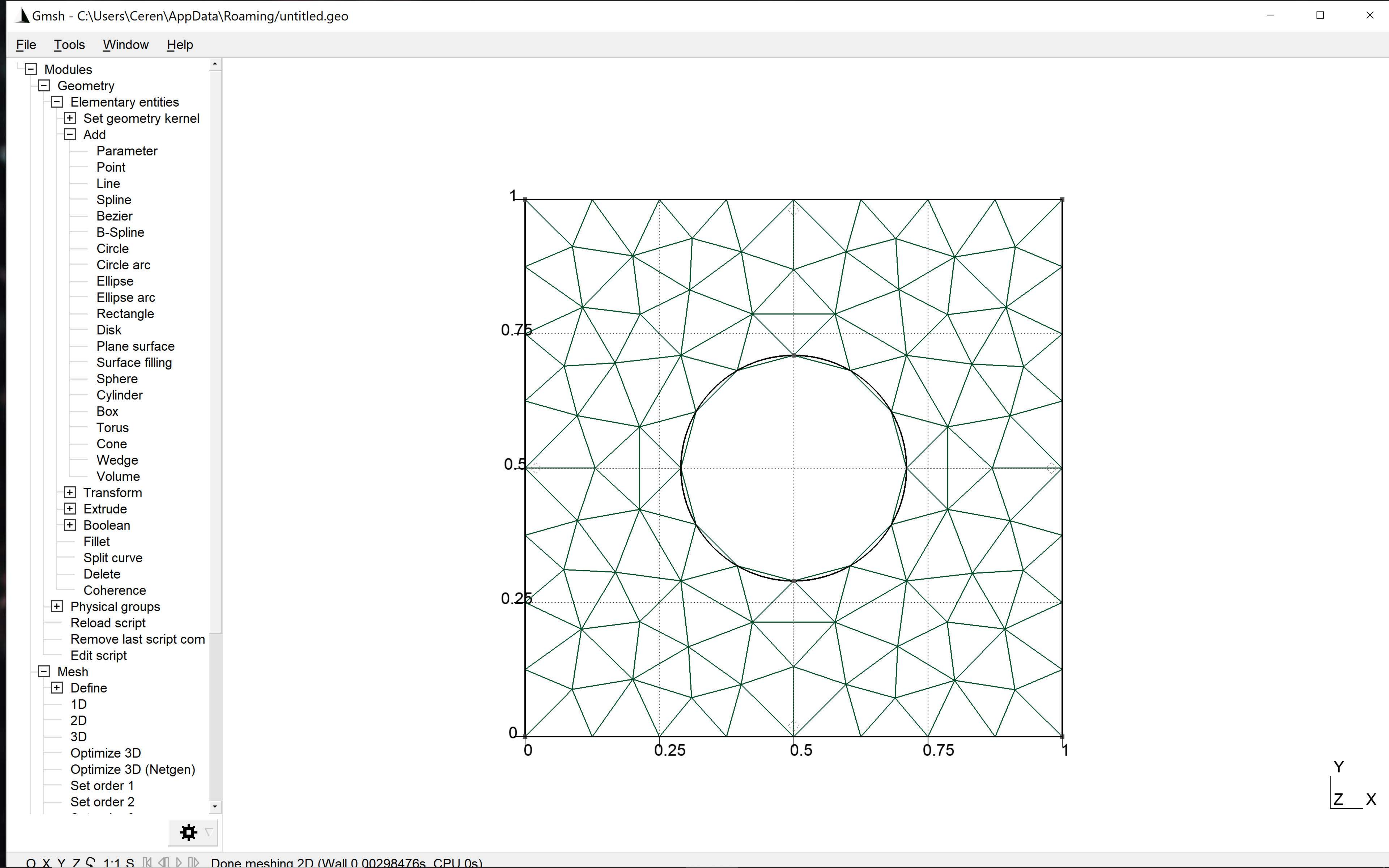}
}
\caption{XHDG mesh with background mesh cut by the boundary (left). Fitted mesh with inner circular boundary necessary for any standard discretization (right).}
\label{fig:DirichletVoidExample}
\end{figure}

\begin{figure}[h]
\centering
\subfigure[]{
\includegraphics[width=0.45\textwidth]{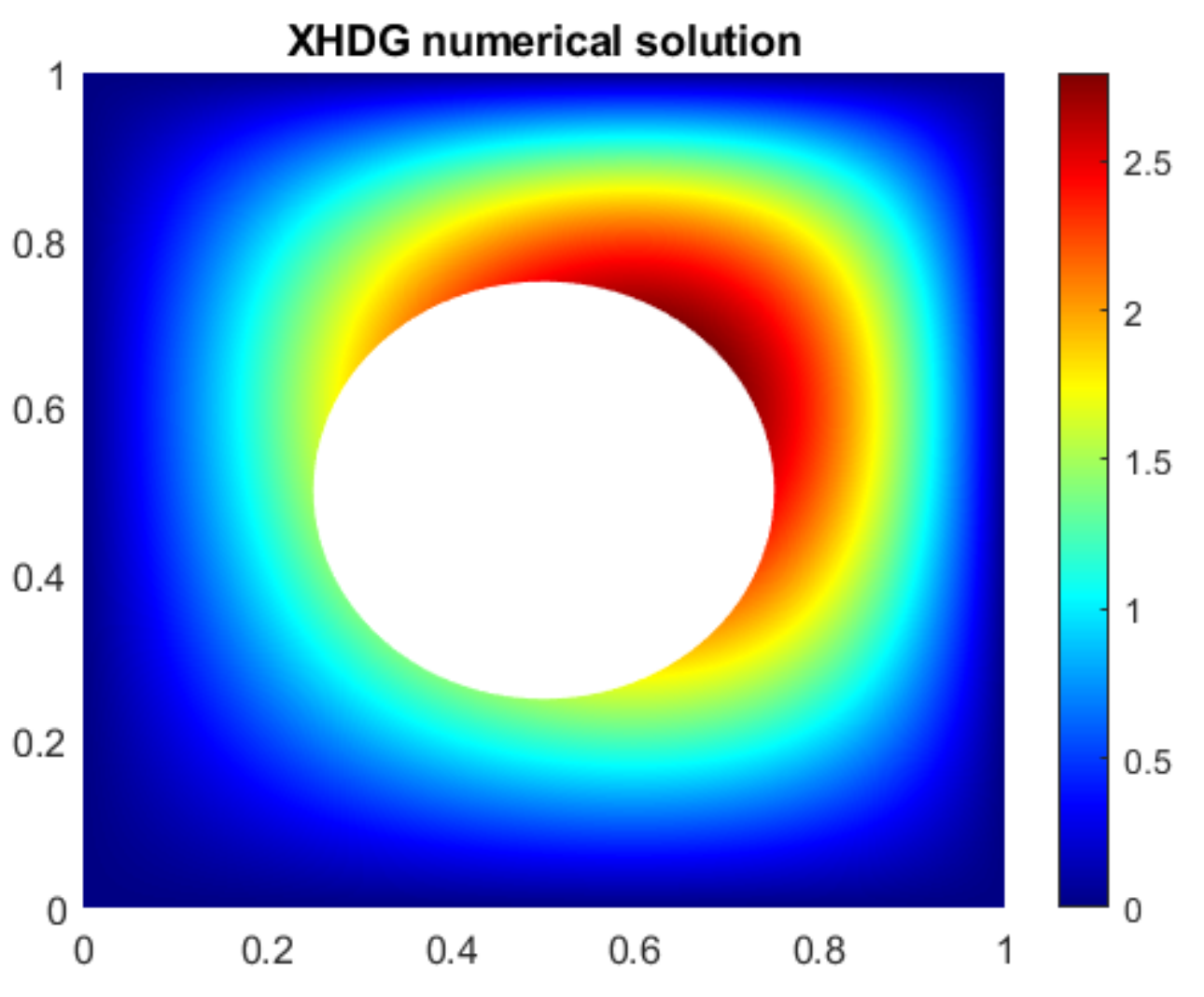}
} 
\subfigure[]{
\includegraphics[width=0.46\textwidth]{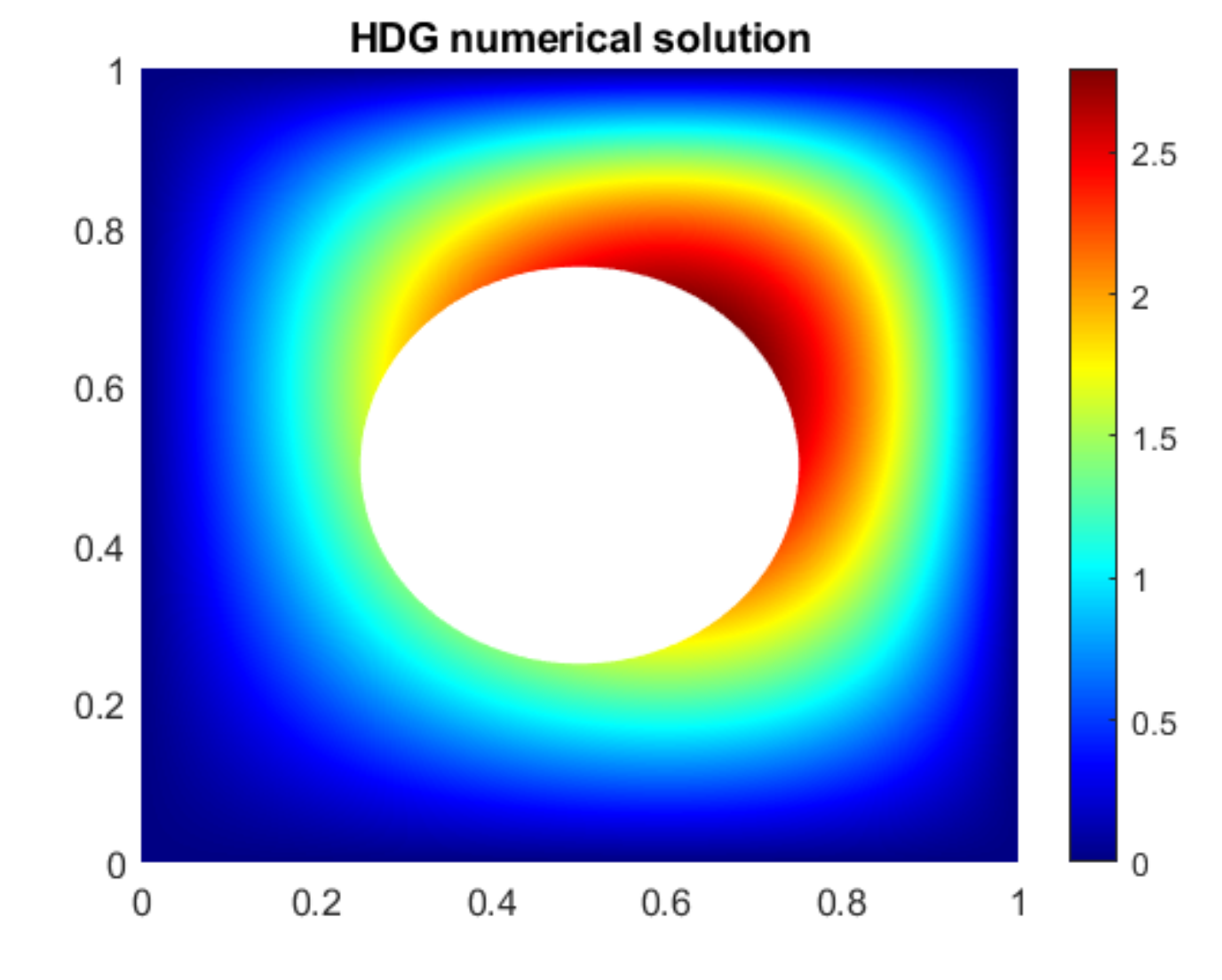  }
}
\caption{XHDG solution over circular cut domain (left), HDG solution over fitted counterpart (right), over diffusion dominated flow regime.}
\label{fig:DirichletVoidSolnDiff}
\end{figure}

For the first numerical example, CD equation defined in \eqref{eq:cd_eqn} is solved over square background mesh with a circular inner void boundary $\Omega=(0,1)^2\backslash B((0.5,0.5),0.42)$. Dirichlet boundary conditions are imposed at the outer and inner cut boundary. The viscosity set to be $\nu=1$ and the convective velocity $c=(1,1)$; leading to a diffusion dominated flow regime. The source term $f$ is set so that the analytical solution is:
$u=\exp(x+y)\sin(\pi x)\sin(\pi y).$
Figure \ref{fig:DirichletVoidExample} shows the XHDG mesh with circular cut boundary and its mesh fitted counterpart. In Figure \ref{fig:DirichletVoidSolnDiff}, the XHDG numerical solution over cut domain (left panel) and HDG solution over its mesh-fitted counterpart is presented. 
%
In Tables \ref{Table_XHDG_VOID_Dirchlet_central_diffusion_diffdom} and \ref{Table_XHDG_VOID_Dirchlet_upwind_diffusion_diffdom} the convergence results of the example is summarized. In Table \ref{Table_XHDG_VOID_Dirchlet_central_diffusion_diffdom} results when centered flux definition; in Table \ref{Table_XHDG_VOID_Dirchlet_upwind_diffusion_diffdom}  results when upwind flux definition used are presented. It can be seen clearly that optimal and super convergence is reached for approximation order $p\leq4$ using XHDG discretization over the cut geometry.  
%
%
\begin{table}[H]
\begin{center}
\caption[Dirichlet circular void example ($D=0.42$): convergence of the X-HDG centered scheme for the steady diffusion-dominated problem]{Dirichlet circular void example ($D=0.42$): convergence of the X-HDG centered scheme for the steady diffusion-dominated problem}\label{Table_XHDG_VOID_Dirchlet_central_diffusion_diffdom}
\resizebox{0.81\textwidth}{!}{%
\begin{tabular}{ccccccc}
\hline
\hline
Degree&Mesh&\multicolumn{2}{c}{$||u-u_{h}||_{L^{2}(\Omega)}$}&&\multicolumn{2}{c}{$||u-u^{\ast}_{h}||_{L^{2}(\Omega)}$}\\ 
\hline
(p)&(n)&error&order&&error&order\\
\hline
$1$&$4$&$8.14e-02$&$-$&&$1.33e-02$&$-$\\
   &$8$&$2.23e-02$&$1.87$&&$1.77e-03$&$2.91$\\
   &$16$&$5.79e-03$&$1.94$&&$2.33e-04$&$2.92$\\
   &$32$&$1.47e-03$&$1.98$&&$2.94e-05$&$2.98$\\
   &$64$&$3.71e-04$&$1.99$&&$3.71e-06$&$2.99$\\
\hline
$2$&$4$&$1.15e-02$&$-$&&$1.25e-03$&$-$\\
   &$8$&$1.5e-03$&$2.93$&&$7.7e-05$&$4.03$\\
   &$16$&$1.93e-04$&$2.96$&&$4.86e-06$&$3.98$\\
   &$32$&$2.45e-05$&$2.98$&&$3.06e-07$&$3.99$\\
   &$64$&$3.08e-06$&$2.99$&&$1.92e-08$&$4.00$\\
\hline
$3$&$4$&$8.58e-04$&$-$&&$8.14e-05$&$-$\\
   &$8$&$5.71e-05$&$3.91$&&$2.58e-06$&$4.98$\\
   &$16$&$3.69e-06$&$3.95$&&$8.22e-08$&$4.97$\\
   &$32$&$2.35e-07$&$3.97$&&$2.61e-09$&$4.98$\\
   &$64$&$1.49e-08$&$3.98$&&$1.1e-09$&$-$\\
\hline
$4$&$4$&$7.69e-05$&$-$&&$6.06e-06$&$-$\\
   &$8$&$2.5e-06$&$4.94$&&$9.52e-08$&$5.99$\\
   &$16$&$8.04e-08$&$4.96$&&$1.52e-09$&$5.97$\\
   &$32$&$2.54e-09$&$4.98$&&$2.39e-11$&$5.99$\\
   &$64$&$1.44e-06$&$-$&&$1.61e-08$&$-$\\
\hline
\hline
\end{tabular}}
\end{center}
\end{table}

\begin{table}[H]
\begin{center}
\caption[Dirichlet circular void example ($D=0.42$): convergence of the X-HDG upwind scheme for the steady diffusion-dominated problem]{Dirichlet circular void example ($D=0.42$): convergence of the X-HDG upwind scheme for the steady diffusion-dominated problem}\label{Table_XHDG_VOID_Dirchlet_upwind_diffusion_diffdom}
\resizebox{0.81\textwidth}{!}{%
\begin{tabular}{ccccccc}
\hline
\hline
Degree&Mesh&\multicolumn{2}{c}{$||u-u_{h}||_{L^{2}(\Omega)}$}&&\multicolumn{2}{c}{$||u-u^{\ast}_{h}||_{L^{2}(\Omega)}$}\\ 
\hline
(p)&(n)&error&order&&error&order\\
\hline
$1$&$4$&$1.39e-01$&$-$&&$1.33e-02$&$-$\\
   &$8$&$3.7e-02$&$1.91$&&$1.75e-03$&$2.92$\\
   &$16$&$9.52e-03$&$1.96$&&$2.28e-04$&$2.94$\\
   &$32$&$2.4e-03$&$1.99$&&$2.8e-05$&$3.02$\\
   &$64$&$6.03e-04$&$1.99$&&$3.51e-06$&$2.99$\\
\hline
$2$&$4$&$1.47e-02$&$-$&&$1.21e-03$&$-$\\
   &$8$&$2.1e-03$&$2.8$&&$7.67e-05$&$3.98$\\
   &$16$&$2.75e-04$&$2.94$&&$4.86e-06$&$3.98$\\
   &$32$&$3.46e-05$&$2.99$&&$3.06e-07$&$3.99$\\
   &$64$&$4.38e-06$&$2.98$&&$1.92e-08$&$4.00$\\
\hline
$3$&$4$&$1.44e-03$&$-$&&$7.4e-05$&$-$\\
   &$8$&$1.01e-04$&$3.82$&&$2.5e-06$&$4.89$\\
   &$16$&$6.61e-06$&$3.94$&&$8.08e-08$&$4.95$\\
   &$32$&$4.17e-07$&$3.99$&&$2.56e-09$&$4.98$\\
   &$64$&$2.82e-05$&$-$&&$2.53e-09$&$-$\\
\hline
$4$&$4$&$1.25e-04$&$-$&&$6.57e-06$&$-$\\
   &$8$&$4.05e-06$&$4.94$&&$1.04e-07$&$5.98$\\
   &$16$&$1.32e-07$&$4.94$&&$1.66e-09$&$5.97$\\
   &$32$&$4.27e-09$&$4.95$&&$2.63e-11$&$5.98$\\
   &$64$&$2.81e-02$&$-$&&$8.46e-08$&$-$\\
\hline
\hline
\end{tabular}}
\end{center}
\end{table}

Additionally here; Tables \ref{Table_XHDG_VOID_Neumann_central_diffusion_diffdom} and \ref{Table_XHDG_VOID_Neumann_upwind_diffusion_diffdom} are added; showing the results when Neumann boundary conditions \eqref{eq:cd_bc} are imposed at the inner cut boundary $\I$. Results show that independent of the boundary condition, optimal and super convergence is ensured using unfitted XHDG discretization. 

\begin{table}[H]
\begin{center}
\caption[Neumann circular void example ($D=0.42$): convergence of the X-HDG centered scheme for the steady diffusion-dominated problem]{Neumann circular void example ($D=0.42$): convergence of the X-HDG centered scheme for the steady diffusion-dominated problem}\label{Table_XHDG_VOID_Neumann_central_diffusion_diffdom}
\resizebox{0.81\textwidth}{!}{%
\begin{tabular}{ccccccc}
\hline
\hline
Degree&Mesh&\multicolumn{2}{c}{$||u-u_{h}||_{L^{2}(\Omega)}$}&&\multicolumn{2}{c}{$||u-u^{\ast}_{h}||_{L^{2}(\Omega)}$}\\ 
\hline
(p)&(n)&error&order&&error&order\\
\hline
$1$&$4$&$8.10e-02$&$-$&&$1.35e-02$&$-$\\
   &$8$&$2.23e-02$&$1.86$&&$1.79e-03$&$2.92$\\
   &$16$&$5.79e-03$&$1.94$&&$2.33e-04$&$2.93$\\
   &$32$&$1.47e-03$&$1.98$&&$2.94e-05$&$3.00$\\
   &$64$&$3.71e-04$&$1.99$&&$3.68e-06$&$2.99$\\
\hline
$2$&$4$&$1.14e-02$&$-$&&$1.25e-03$&$-$\\
   &$8$&$1.50e-03$&$2.93$&&$7.77e-05$&$4.02$\\
   &$16$&$1.93e-04$&$2.96$&&$4.89e-06$&$3.99$\\
   &$32$&$2.45e-05$&$2.98$&&$3.08e-07$&$3.99$\\
   &$64$&$3.08e-06$&$2.99$&&$1.92e-08$&$3.99$\\
\hline
$3$&$4$&$8.58e-04$&$-$&&$8.14e-05$&$-$\\
   &$8$&$5.71e-05$&$3.91$&&$2.59e-06$&$4.98$\\
   &$16$&$3.69e-06$&$3.95$&&$8.23e-08$&$4.98$\\
   &$32$&$2.35e-07$&$3.97$&&$2.61e-09$&$4.98$\\
   &$64$&$1.48e-08$&$3.99$&&$7.61e-10$&$-$\\
\hline
$4$&$4$&$7.69e-05$&$-$&&$6.08e-06$&$-$\\
   &$8$&$2.50e-06$&$4.94$&&$9.52e-08$&$5.99$\\
   &$16$&$8.04e-08$&$4.96$&&$1.52e-09$&$5.97$\\
   &$32$&$2.54e-09$&$4.98$&&$2.39e-11$&$5.98$\\
   &$64$&$8.06e-06$&$-$&&$9.01e-09$&$-$\\
\hline
\hline
\end{tabular}}
\end{center}
\end{table}

\begin{table}[H]
\begin{center}
\caption[Neumann circular void example ($D=0.42$): convergence of the X-HDG upwind scheme for the steady diffusion-dominated problem]{Neumann circular void example ($D=0.42$): convergence of the X-HDG upwind scheme for the steady diffusion-dominated problem}\label{Table_XHDG_VOID_Neumann_upwind_diffusion_diffdom}
\resizebox{0.81\textwidth}{!}{%
\begin{tabular}{ccccccc}
\hline
\hline
Degree&Mesh&\multicolumn{2}{c}{$||u-u_{h}||_{L^{2}(\Omega)}$}&&\multicolumn{2}{c}{$||u-u^{\ast}_{h}||_{L^{2}(\Omega)}$}\\ 
\hline
(p)&(n)&error&order&&error&order\\
\hline
$1$&$4$&$1.37e-01$&$-$&&$1.37e-02$&$-$\\
   &$8$&$3.69e-02$&$1.89$&&$1.87e-03$&$2.87$\\
   &$16$&$9.51e-03$&$1.96$&&$2.45e-04$&$2.93$\\
   &$32$&$2.40e-03$&$1.99$&&$2.98e-05$&$3.04$\\
   &$64$&$6.03e-04$&$1.99$&&$3.69e-06$&$3.01$\\
\hline
$2$&$4$&$1.59e-02$&$-$&&$1.26e-03$&$-$\\
   &$8$&$2.09e-03$&$2.93$&&$7.81e-05$&$4.02$\\
   &$16$&$2.77e-04$&$2.91$&&$4.93e-06$&$3.98$\\
   &$32$&$3.59e-05$&$2.95$&&$3.12e-07$&$3.98$\\
   &$64$&$4.52e-06$&$2.99$&&$1.95e-08$&$4.00$\\
\hline
$3$&$4$&$1.44e-03$&$-$&&$7.43e-05$&$-$\\
   &$8$&$9.92e-05$&$3.86$&&$2.50e-06$&$4.89$\\
   &$16$&$3.72e-06$&$3.88$&&$8.11e-08$&$4.95$\\
   &$32$&$4.23e-07$&$3.99$&&$2.57e-09$&$4.98$\\
   &$64$&$1.75e-04$&$-$&&$6.04e-09$&$-$\\
\hline
$4$&$4$&$1.27e-04$&$-$&&$6.64e-06$&$-$\\
   &$8$&$4.23e-06$&$4.91$&&$1.05e-07$&$5.98$\\
   &$16$&$1.38e-07$&$4.94$&&$1.67e-09$&$5.97$\\
   &$32$&$4.87e-09$&$4.83$&&$2.84e-11$&$5.88$\\
   &$64$&$3.49e-01$&$-$&&$7.64e-07$&$-$\\
\hline
\hline
\end{tabular}}
\end{center}
\end{table}

\begin{remark}\label{rm:BadCut}
In Tables \ref{Table_XHDG_VOID_Dirchlet_central_diffusion_diffdom},
\ref{Table_XHDG_VOID_Dirchlet_upwind_diffusion_diffdom},
\ref{Table_XHDG_VOID_Neumann_central_diffusion_diffdom} and \ref{Table_XHDG_VOID_Neumann_upwind_diffusion_diffdom} the reader can notice that for a high approximation degree and fine mesh, although the optimal and super convergence is reached, there might be a jump in error value and hence, the convergence order is not shown. This is the main bottleneck of any unfitted discretization, namely, the ill conditioning because of so called \textit{bad-cut} situations. \textit{Bad-cut} situation occurs when the boundary cuts the element in the background mesh very close to its vertex or edge, leading to a very small contribution from that element in the global system. At the end, these small contributions lead to a global system matrix with a high condition number, effecting the error and convergence of the cut discretization scheme. It is not in the context of this work to propose a remedy to this ill conditioning problem, however, the reader is referred to \cite{Gurkan-Kronbichler-Fernandez:17} to a node-moving type quick-fix or to \cite{Gurkan-Sticko-Massing:20} to a ghost penalty type solution to this ill conditioning problem.      
\end{remark}
\subsection{Circular cut boundary on a square domain - Convection dominated}\label{sec:SecondExample}
\begin{figure}[h]
\centering
\subfigure[]{
\includegraphics[width=0.46\textwidth]{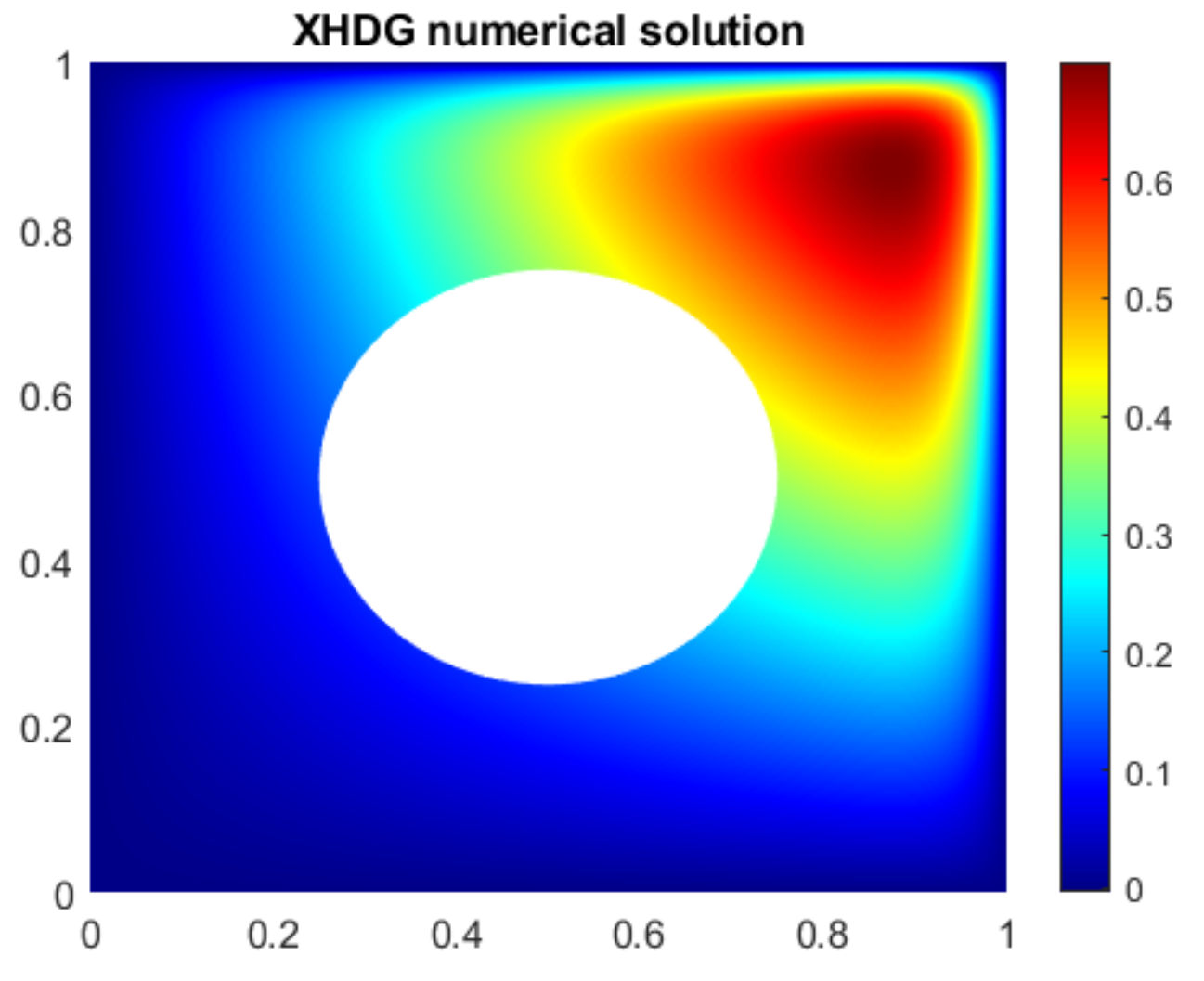}
}
\subfigure[]{
\includegraphics[width=0.47\textwidth]{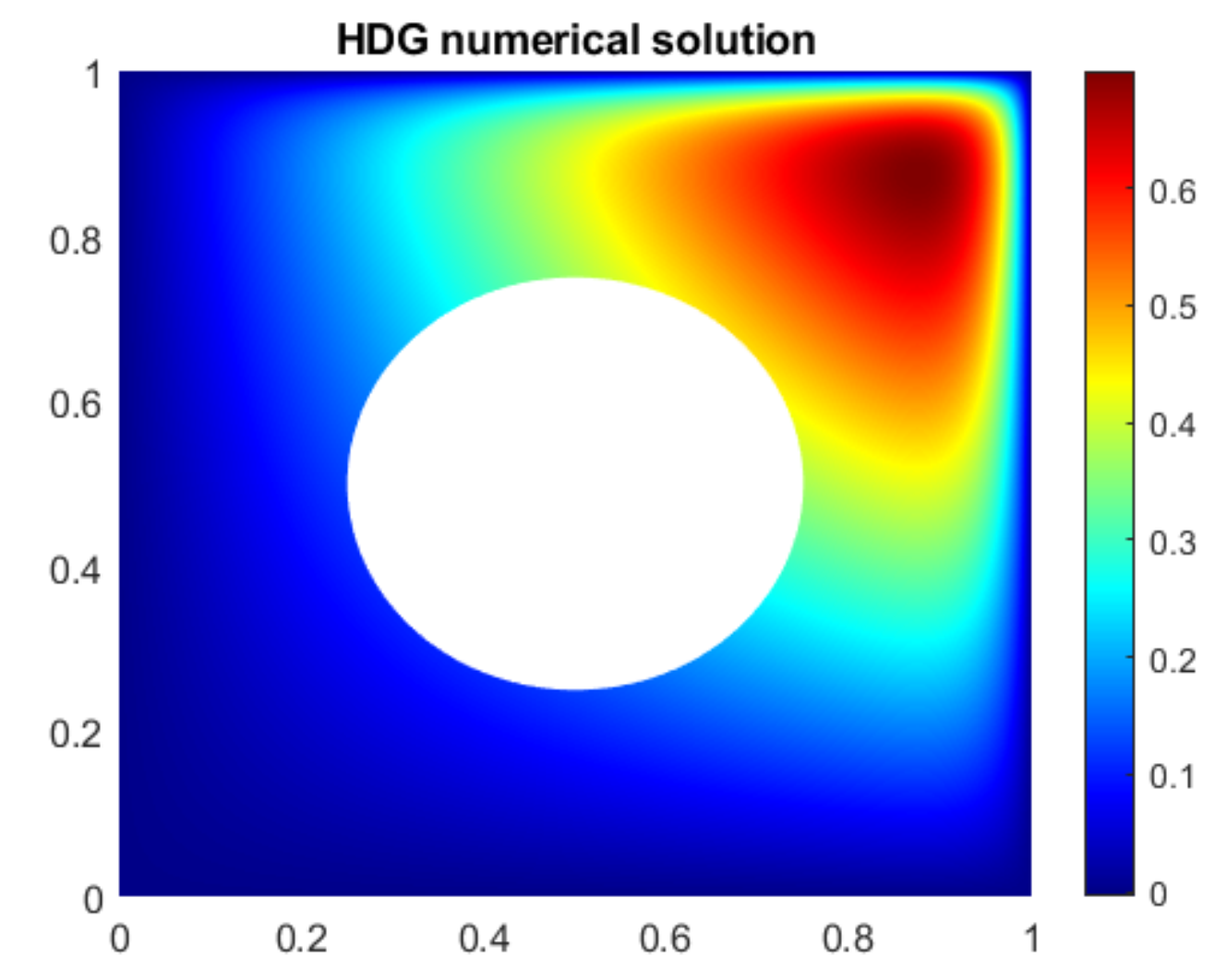}
}
\caption{XHDG solution over circular cut domain (left), HDG solution over fitted counterpart (right), over convection dominated flow regime.}
\label{fig:DirichletVoidSolnConv}
\end{figure}

In this second numerical test, the domain and the cut boundary kept the same as in first example, as shown in Figure \ref{fig:DomainVithCircularVoidExample} left panel. The diffusion coefficient reduced to $\nu=0.05$ while keeping the convective velocity as $c=(1,1)$; leading to a convection dominated flow regime with $Pe=20$ where $Pe:=\frac{cL}{\nu}$ and $L$ is the edge length of the square domain. Dirichlet boundary conditions are imposed at the cut boundary and the analytical solution is set to be
\begin{equation*}
u(x,y)=xy\frac{(1-e^{(x-1)c_x}(1-e^{(y-1)c_y})}{(1-e^{-c_x})(1-e^{-c_y})}
\end{equation*} 
For high $Pe$ number this analytical solution creates boundary layers at $x=1$ and $y=1$. The XHDG solution over the cut domain and HDG solution over its fitted counterpart for this convection dominated flow is shown in Figure \ref{fig:DirichletVoidSolnConv}.
Tables \ref{Table_lowVISC_XHDG_VOID_Dirichlet_central_convection_convdom} and \ref{Table_lowVISC_XHDG_VOID_Dirichlet_upwind_convection_convdom} shows the convergence and accuracy of XHDG for convection dominated flow regime when centered and upwind flux is used, respectively. As the results show, XHDG reaches optimal and super convergence over a convection dominated flow regime when boundary layers exists in the domain as well. 
%
\begin{table}[H]
\begin{center}
\caption{Dirichlet circular void example ($D=0.42$, $\nu=0.05$ and $c=(1,1)$): convergence of the X-HDG centered scheme for the steady convection-dominated problem}\label{Table_lowVISC_XHDG_VOID_Dirichlet_central_convection_convdom}
\resizebox{0.81\textwidth}{!}{%
\begin{tabular}{ccccccc}
\hline
\hline
Degree&Mesh&\multicolumn{2}{c}{$||u-u_{h}||_{L^{2}(\Omega)}$}&&\multicolumn{2}{c}{$||u-u^{\ast}_{h}||_{L^{2}(\Omega)}$}\\ 
\hline
(p)&(n)&error&order&&error&order\\
\hline
$1$&$4$&$1.77e-03$&$-$&&$9.59e-04$&$-$\\
   &$8$&$4.46e-04$&$1.99$&&$1.85e-04$&$2.4$\\
   &$16$&$1.10e-04$&$2.01$&&$3.15e-05$&$2.6$\\
   &$32$&$2.73e-05$&$2.01$&&$4.62e-06$&$2.8$\\ 
   &$64$&$6.80e-06$&$2.03$&&$6.38e-07$&$2.9$\\     
\hline
$2$&$4$&$2.92e-04$&$-$&&$1.42e-04$&$-$\\
   &$8$&$3.47e-05$&$3.07$&&$1.09e-05$&$3.7$\\
   &$16$&$4.28e-06$&$3.02$&&$7.71e-07$&$3.8$\\
   &$32$&$5.34e-07$&$3.01$&&$5.13e-08$&$3.9$\\ 
   &$64$&$6.69e-08$&$3.00$&&$3.31e-09$&$3.9$\\     
\hline
$3$&$4$&$2.07e-05$&$-$&&$8.49e-06$&$-$\\
   &$8$&$1.24e-06$&$4.07$&&$3.33e-07$&$3.7$\\
   &$16$&$7.65e-08$&$4.01$&&$1.18e-08$&$3.8$\\
   &$32$&$4.79e-09$&$4.00$&&$3.93e-10$&$3.9$\\ 
   &$64$&$3.04e-10$&$3.98$&&$2.79e-11$&$3.9$\\
\hline
$4$&$4$&$7.76e-07$&$-$&&$2.79e-07$&$-$\\
   &$8$&$2.34e-08$&$5.05$&&$5.56e-09$&$5.7$\\
   &$16$&$7.30e-10$&$5.00$&&$9.91e-11$&$5.8$\\
   &$32$&$2.29e-11$&$4.99$&&$1.65e-12$&$5.9$\\ 
   &$64$&$1.62e-06$&$-$&&$1.28e-09$&$-$\\
\hline
\hline
\end{tabular}}
\end{center}
\end{table}

\begin{table}[H]
\begin{center}
\caption{Dirichlet circular void example ($D=0.42$, $\nu=0.05$ and $c=(1,1)$): convergence of the X-HDG upwind scheme for the steady convection-dominated problem}\label{Table_lowVISC_XHDG_VOID_Dirichlet_upwind_convection_convdom}
\resizebox{0.81\textwidth}{!}{%
\begin{tabular}{ccccccc}
\hline
\hline
Degree&Mesh&\multicolumn{2}{c}{$||u-u_{h}||_{L^{2}(\Omega)}$}&&\multicolumn{2}{c}{$||u-u^{\ast}_{h}||_{L^{2}(\Omega)}$}\\ 
\hline
(p)&(n)&error&order&&error&order\\
\hline
$1$&$4$&$1.96e-03$&$-$&&$5.41e-04$&$-$\\
   &$8$&$5.22e-04$&$1.91$&&$9.34e-05$&$2.5$\\
   &$16$&$1.33e-04$&$1.96$&&$1.48e-05$&$2.7$\\
   &$32$&$3.38e-05$&$1.98$&&$2.09e-06$&$2.8$\\ 
   &$64$&$8.47e-06$&$1.99$&&$2.77e-07$&$2.9$\\     
\hline
$2$&$4$&$3.27e-04$&$-$&&$8.11e-05$&$-$\\
   &$8$&$4.12e-05$&$2.99$&&$5.84e-06$&$3.79$\\
   &$16$&$5.15e-06$&$3.0$&&$4.01e-07$&$3.87$\\
   &$32$&$6.43e-07$&$3.0$&&$2.62e-08$&$3.94$\\ 
   &$64$&$8.04e-08$&$3.0$&&$1.67e-09$&$3.97$\\     
\hline
$3$&$4$&$2.48e-05$&$-$&&$3.53e-06$&$-$\\
   &$8$&$1.56e-06$&$3.99$&&$1.21e-07$&$4.87$\\
   &$16$&$9.78e-08$&$3.99$&&$4.06e-09$&$4.89$\\
   &$32$&$6.13e-09$&$4.00$&&$1.32e-10$&$4.95$\\ 
   &$64$&$4.77e-08$&$-$&&$2.38e-10$&$-$\\
\hline
$4$&$4$&$9.88e-07$&$-$&&$5.38e-08$&$-$\\
   &$8$&$3.11e-08$&$4.99$&&$3.68e-10$&$5.95$\\
   &$16$&$9.80e-10$&$4.99$&&$1.49e-11$&$5.87$\\
   &$32$&$3.08e-11$&$4.99$&&$2.26e-13$&$6.04$\\ 
   &$64$&$9.06e-04$&$-$&&$1.43e-09$&$-$\\
\hline
\hline
\end{tabular}}
\end{center}
\end{table}

\subsection{Peanut shaped cut boundary on a square domain - Convection dominated}
\begin{figure}[h]
\centering
\subfigure[]{
\includegraphics[width=0.40\textwidth]{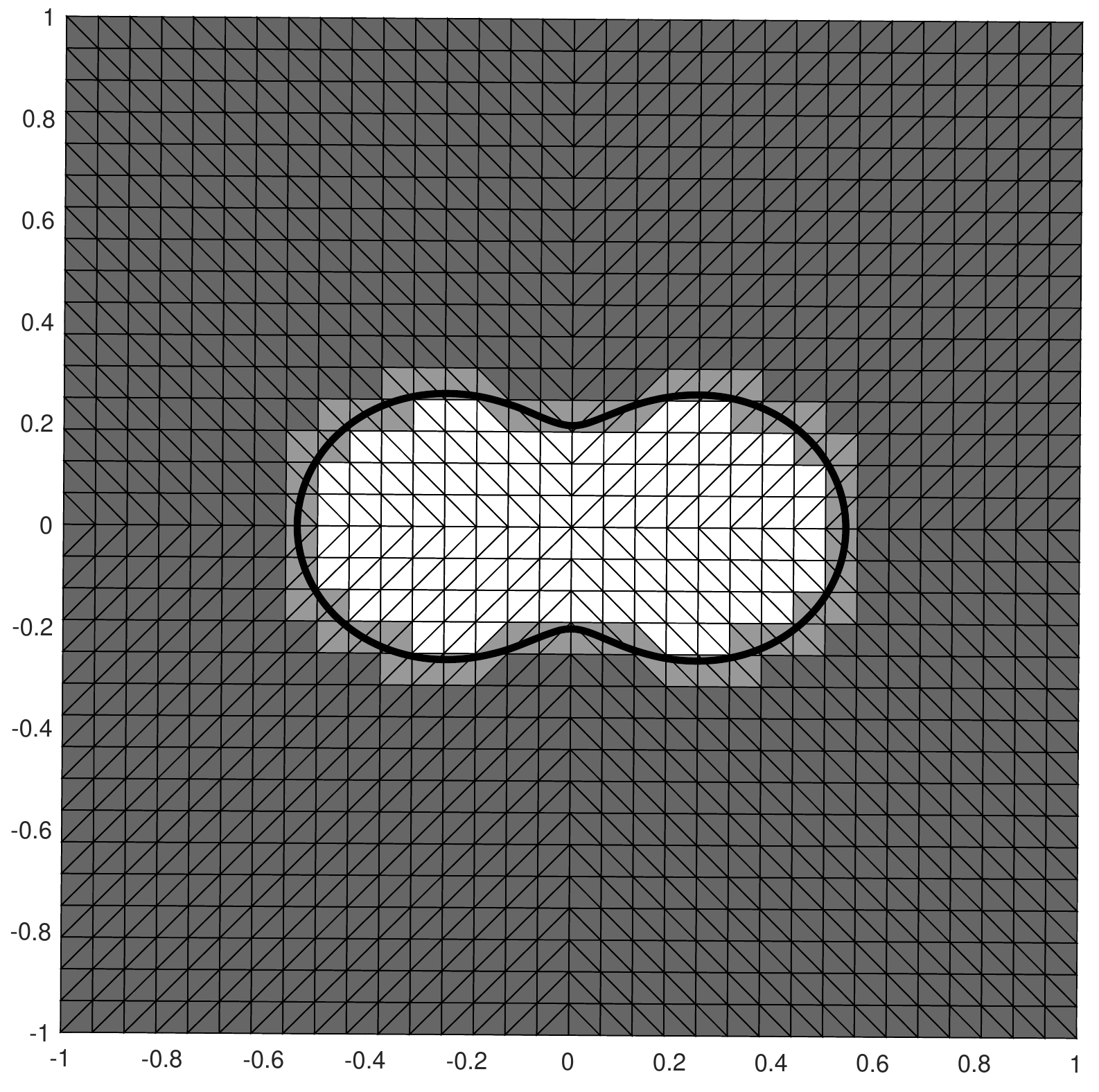}
}
\subfigure[]{
\includegraphics[width=0.53\textwidth]{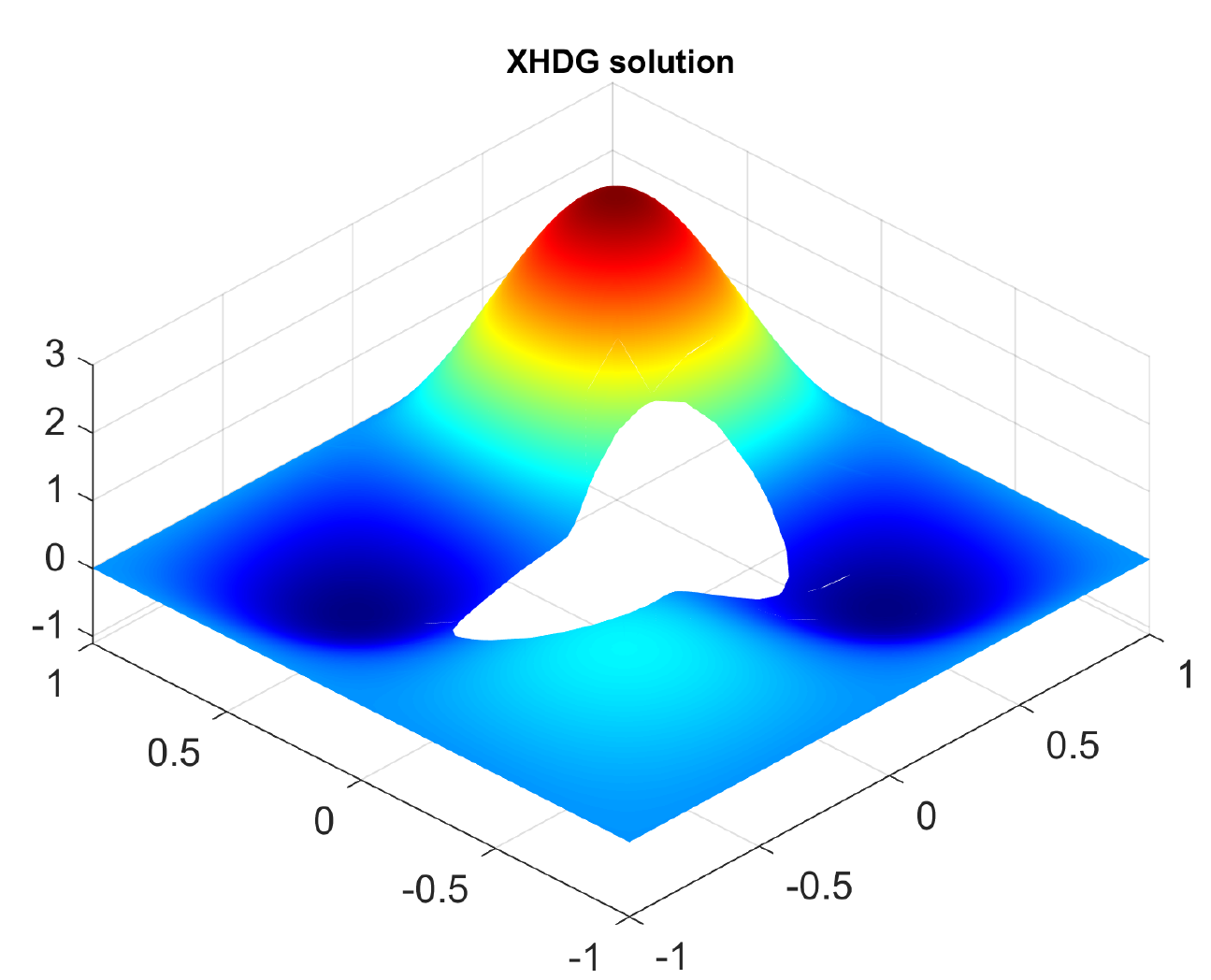}
}
\caption{The square background mesh with peanut cut boundary (left). The XHDG solution (right).}
\label{fig:Peanut_boundary}
\end{figure}

As the last steady example, we have studied a rather challenging peanut shaped cut boundary in a square background domain as shown in Figure \ref{fig:Peanut_boundary}. The convective velocity and the viscosity are set to be $c=(25,25)$ and $\nu=1$, respectively; leading to a convection dominated flow with $Pe=50$, even higher than in previous example. The background geometry is once more in square shape $\Omega=(-1,1)^2 \backslash \I$ and the level set function that defines the peanut shaped boundary is: 
$\Phi=\sqrt{x^2+y^2}-r_0-r_1 \cos⁡(2\arctan_2(x,y))$ where $r_0=0.37$ and $r_1=0.17$. 
The analytical solution is the same as in the example \ref{sec:SecondExample} and Neumann boundary conditions are applied at the cut boundary $\I$.

%
As usual, the results are summarized in Tables \ref{Table_XHDG_Peanut_Neumann_central_convection_convdom} and \ref{Table_XHDG_Peanut_Neumann_upwind_convection_convdom} for centered and upwind flux, respectively. As it is shown in the Tables; although the relatively complex cut geometry in this example evokes the bad-cut situations mentioned in Remark \ref{rm:BadCut}; the XHDG discretization does converge towards the expected orders of convergence with increasing element size.  
\begin{table}[H]
\begin{center}
\caption{Neumann peanut shaped void example ($\nu=1$ \& $c=(25,25)$): convergence of the X-HDG centered scheme for the steady convection-dominated problem}\label{Table_XHDG_Peanut_Neumann_central_convection_convdom}
\resizebox{0.81\textwidth}{!}{%
\begin{tabular}{ccccccc}
\hline
\hline
Degree&Mesh&\multicolumn{2}{c}{$||u-u_{h}||_{L^{2}(\Omega)}$}&&\multicolumn{2}{c}{$||u-u^{\ast}_{h}||_{L^{2}(\Omega)}$}\\ 
\hline
(p)&(n)&error&order&&error&order\\
\hline
$1$&$4$&$7.55e-02$&$-$&&$1.49e-01$&$-$\\
   &$8$&$6.10e-02$&$0.31$&&$5.90e-02$&$1.33$\\
   &$16$&$3.42e-02$&$0.84$&&$1.81e-02$&$1.71$\\
   &$32$&$1.25e-02$&$1.45$&&$3.95e-03$&$2.19$\\ 
   &$64$&$3.61e-03$&$1.79$&&$6.55e-04$&$2.59$\\     
\hline
$2$&$4$&$6.82e-02$&$-$&&$5.24e-02$&$-$\\
  &$8$&$3.29e-02$&$1.05$&&$2.04e-02$&$1.82$\\
   &$16$&$9.44e-03$&$1.80$&&$3.68e-03$&$2.47$\\
   &$32$&$1.78e-03$&$2.41$&&$4.16e-04$&$3.15$\\ 
   &$64$&$2.61e-04$&$2.77$&&$3.32e-05$&$3.65$\\     
\hline
$3$&$4$&$3.73e-02$&$-$&&$3.55e-02$&$-$\\
  &$8$&$1.19e-02$&$1.64$&&$6.45e-03$&$2.46$\\
   &$16$&$2.03e-03$&$2.56$&&$7.17e-04$&$3.17$\\
   &$32$&$2.17e-04$&$3.23$&&$4.42e-05$&$4.02$\\ 
   &$64$&$1.68e-05$&$3.69$&&$1.78e-06$&$4.63$\\
\hline
$4$&$4$&$2.00e-02$&$-$&&$1.55e-02$&$-$\\
   &$8$&$3.85e-03$&$2.38$&&$2.04e-03$&$2.93$\\
   &$16$&$4.28e-04$&$3.17$&&$1.50e-04$&$3.77$\\
   &$32$&$2.68e-05$&$4.00$&&$5.02e-06$&$4.90$\\ 
   &$64$&$1.09e-06$&$4.62$&&$1.03e-07$&$5.61$\\
\hline
\hline
\end{tabular}}
\end{center}
\end{table}

\begin{table}[H]
\begin{center}
\caption{Neumann peanut shaped void example ($\nu=1$ \& $c=(25,25)$): convergence of the X-HDG upwind scheme for the steady convection-dominated problem}\label{Table_XHDG_Peanut_Neumann_upwind_convection_convdom}
\resizebox{0.71\textwidth}{!}{%
\begin{tabular}{ccccccc}
\hline
\hline
Degree&Mesh&\multicolumn{2}{c}{$||u-u_{h}||_{L^{2}(\Omega)}$}&&\multicolumn{2}{c}{$||u-u^{\ast}_{h}||_{L^{2}(\Omega)}$}\\ 
\hline
(p)&(n)&error&order&&error&order\\
\hline
$1$&$4$&$6.54e-02$&$-$&&$1.56e-01$&$-$\\
   &$8$&$5.98e-02$&$0.13$&&$5.62e-02$&$1.47$\\
   &$16$&$3.36e-02$&$0.83$&&$1.66e-02$&$1.76$\\
   &$32$&$1.23e-02$&$1.45$&&$3.44e-03$&$2.27$\\ 
   &$64$&$3.54e-03$&$1.79$&&$5.44e-04$&$2.66$\\     
\hline
$2$&$4$&$6.75e-02$&$-$&&$7.09e-02$&$-$\\
   &$8$&$3.25e-02$&$1.05$&&$1.94e-02$&$1.87$\\
   &$16$&$9.29e-03$&$1.81$&&$3.36e-03$&$2.53$\\
   &$32$&$1.76e-03$&$2.40$&&$3.68e-04$&$3.19$\\ 
   &$64$&$2.58e-04$&$2.77$&&$2.87e-05$&$3.68$\\     
\hline
$3$&$4$&$3.68e-02$&$-$&&$3.44e-02$&$-$\\
   &$8$&$1.18e-02$&$1.64$&&$6.09e-03$&$2.50$\\
   &$16$&$2.01e-03$&$2.55$&&$6.71e-04$&$3.18$\\
   &$32$&$2.14e-04$&$3.23$&&$4.13e-05$&$4.02$\\ 
   &$64$&$1.66e-05$&$3.69$&&$1.66e-06$&$4.64$\\
\hline
$4$&$4$&$1.97e-02$&$-$&&$1.48e-02$&$-$\\
   &$8$&$3.85e-03$&$2.37$&&$1.95e-03$&$2.93$\\
   &$16$&$4.28e-04$&$3.17$&&$1.45e-04$&$3.75$\\
   &$32$&$2.65e-05$&$4.00$&&$4.88e-06$&$4.89$\\ 
   &$64$&$1.09e-06$&$4.62$&&$9.99e-08$&$5.61$\\
\hline
\hline
\end{tabular}}
\end{center}
\end{table}
\subsection{Time dependent CDE - Transport of a Gaussian Pulse}
\begin{figure*}
\centering
\begin{subfigure}
\centering
\includegraphics[width=0.45 \textwidth]{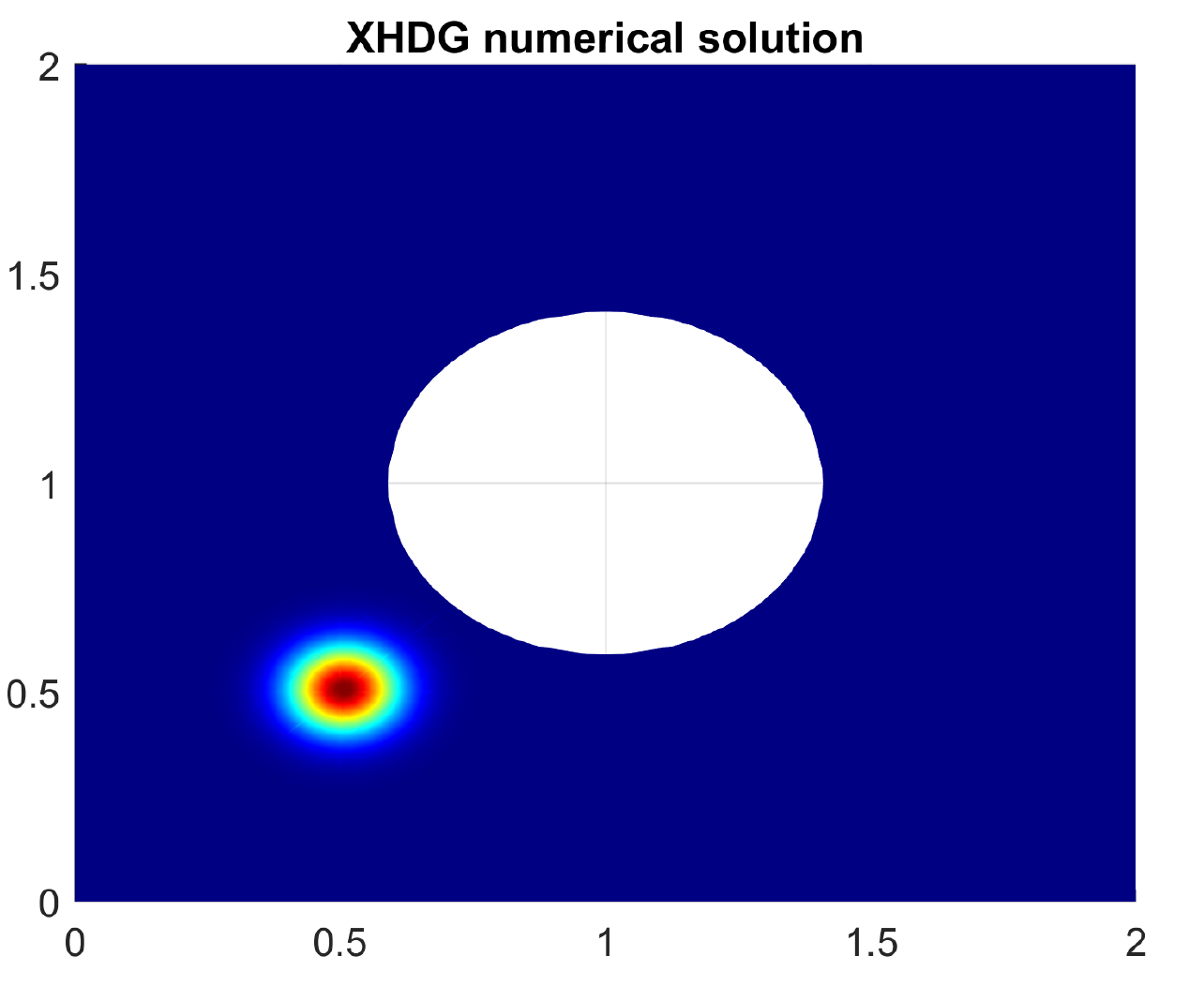}
\end{subfigure}%
\hfill
\begin{subfigure}
\centering
\includegraphics[width=0.47 \textwidth]{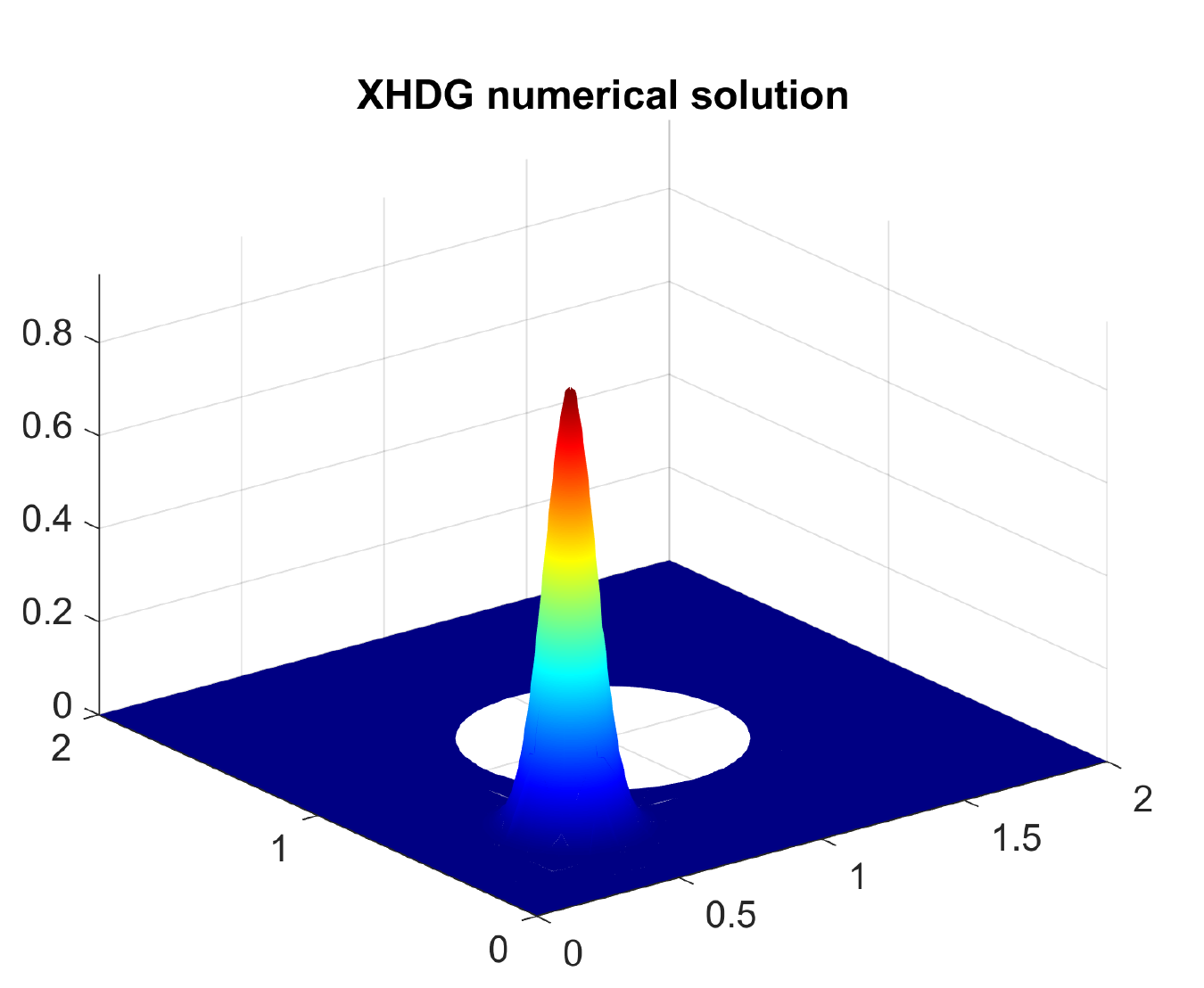}
\end{subfigure}

\bigskip 

\begin{subfigure}
\centering
\includegraphics[width=0.45 \textwidth]{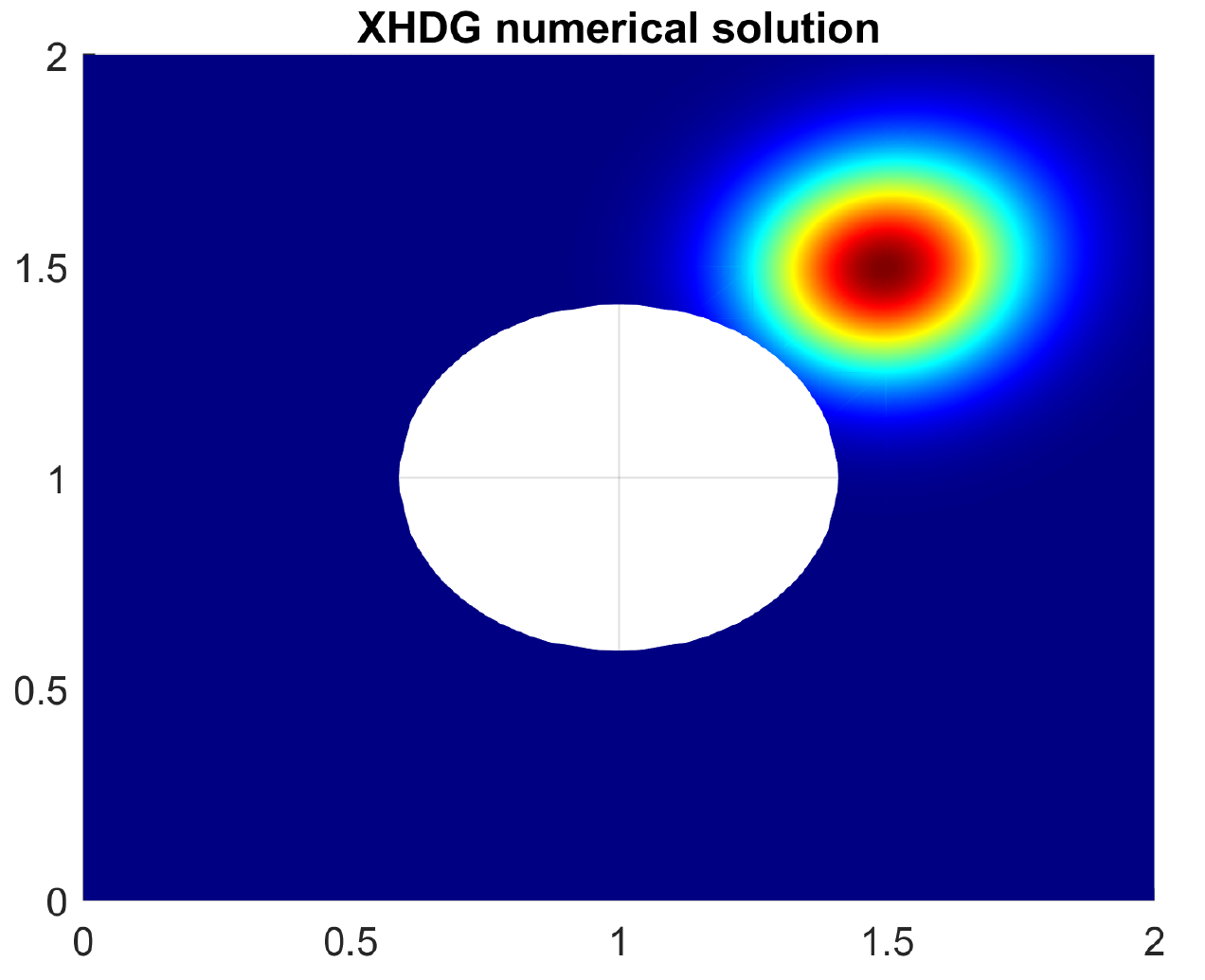}
\end{subfigure}%
\hfill
\begin{subfigure}
\centering
\includegraphics[width=0.47\textwidth]{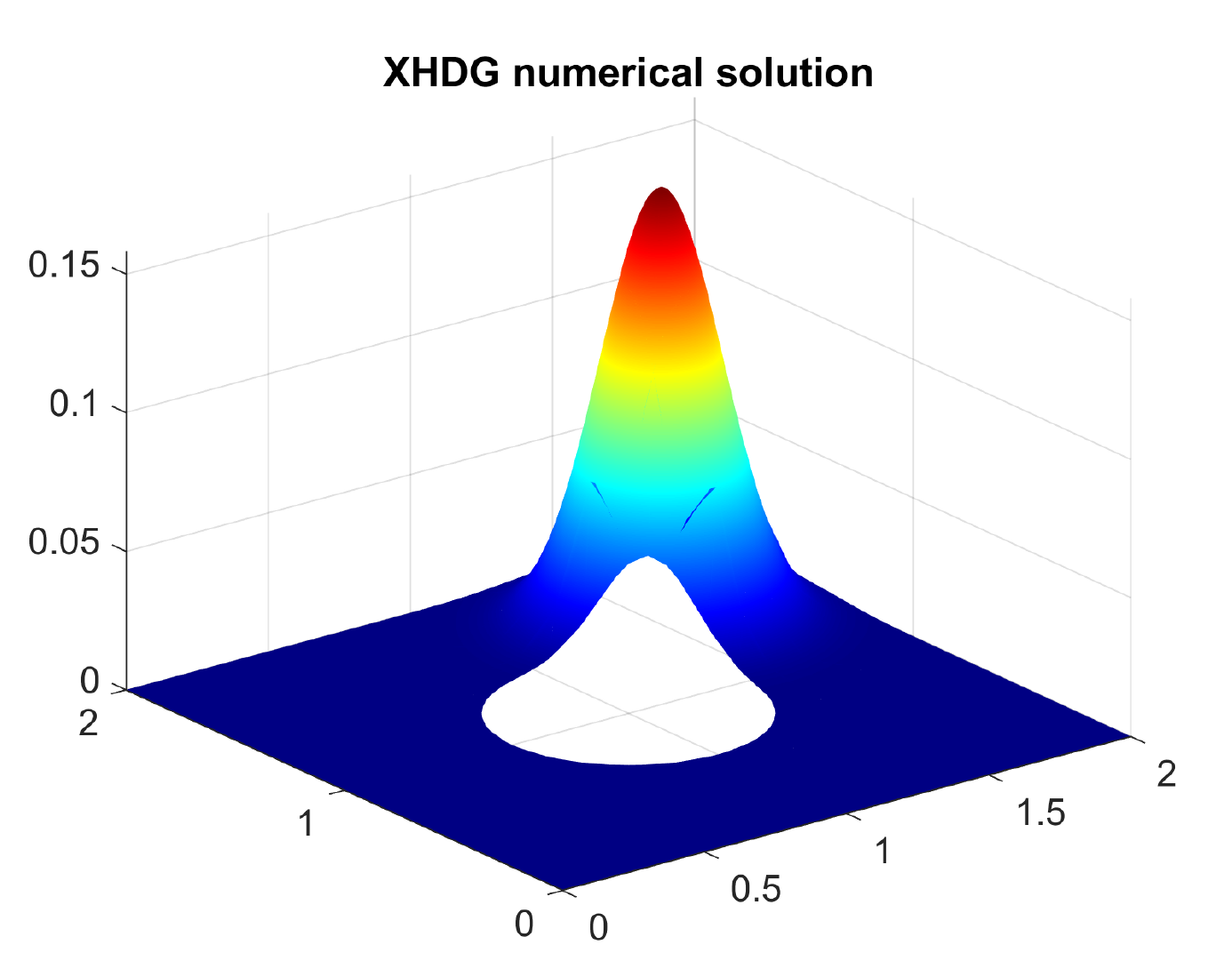}
\end{subfigure}

\caption{The initial position and height (top) and the final $t=1.25$ position and height of Gaussian pulse over cut domain (bottom) obtained using XHDG discretization. The convection and diffusion of pulse can clearly be observed.}
\label{fig:Gaussian_pulse}
\end{figure*}

Last but not least, time dependent CDE is solved over a square background mesh with a circular void $\Omega=(0,2)^2 \backslash B((1,1),0.5)$. Over this domain a Gaussian Pulse is transported following CDE, diagonally. The flow is once more convection dominated with convective velocity $c=(0.8,0.8)$, viscosity $\nu=0.01$ and $Pe=160$.     
The initial condition and the analytical solution are defined as follows:
\begin{equation*}
u(x,y,0)= \exp(-\frac{(x-0.5)^2}{\nu}-\frac{(y-0.5)^2}{\nu})
\end{equation*}
\begin{equation*}
u(x,y,t)= \frac{1}{4t+1} \exp(-\frac{(x-c_x t-0.5)^2}{(4t+1)\nu}-\frac{(y-c_y t-0.5)^2}{(4t+1)\nu})
\end{equation*} 
%
The initial condition leads to a Gaussian pulse of height $h=1$ at position $(x,y)=(0.5,0.5)$ at $t=0$. When $t=1.25$, analytical solution state that the pulse should move to $(x,y)=(1.5,1.5)$ with its height decayed to $h=1/6$. 
Figure \ref{fig:Gaussian_pulse} shows the starting $t=0$ and $t=1.25$ position of the pulse when the problem is solved using XHDG discretization.
 
\begin{figure}[h]
\centering
\includegraphics[width=1.0\textwidth]{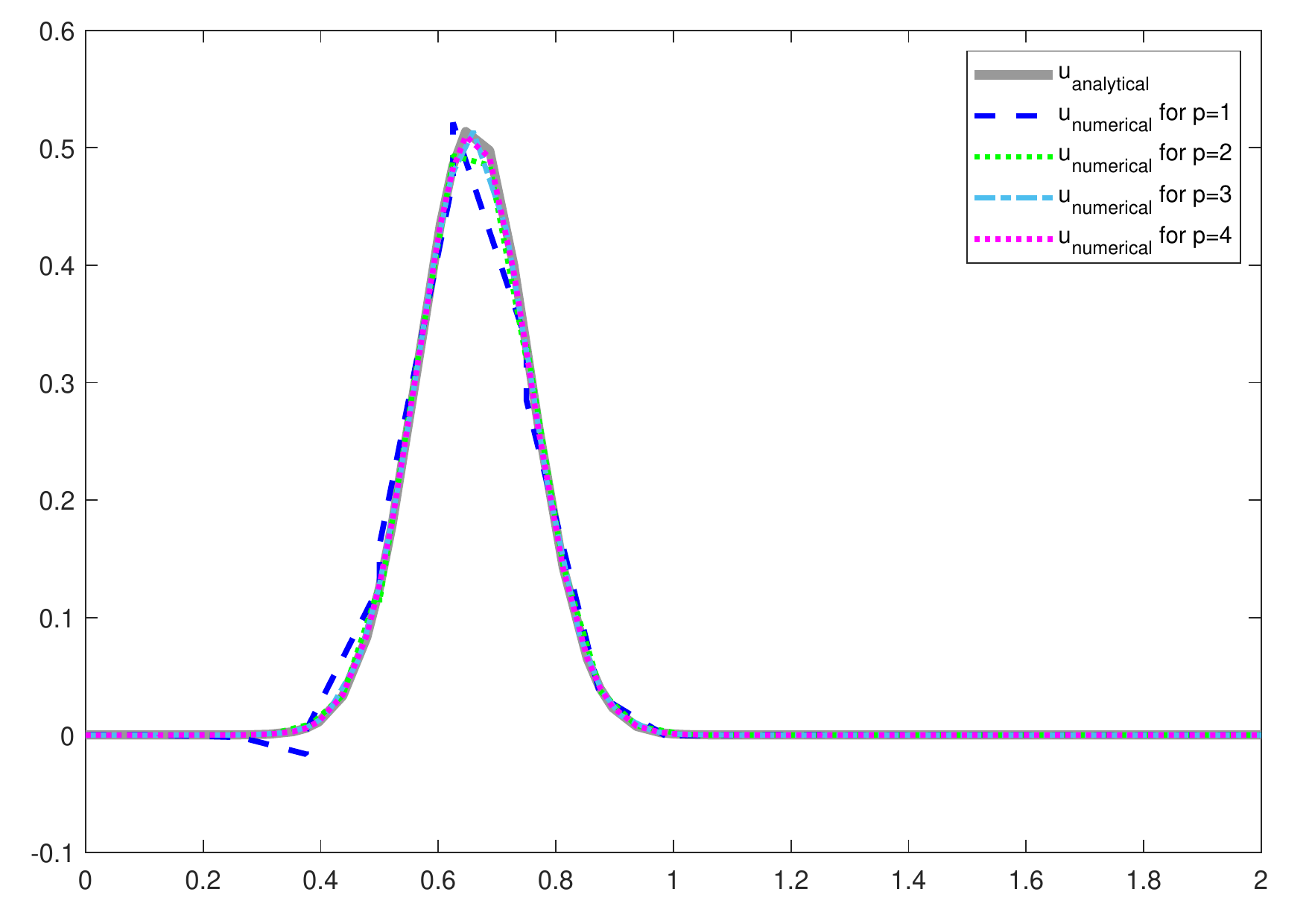}
\caption{The comparison of various XHDG approximations with analytical solution at $x=0.625$ and when $t=0.2$.}
\label{fig:Gauss_comp}
\end{figure}

%
The XHDG solution obtained using $p=1,2,3 \text{ and } 4^{th}$ order approximation functions along $x=0.625$ when $t=0.2$, before pulse entering the void, is compared with analytical solution in Figure \ref{fig:Gauss_comp}. Results show that XHDG solution is in perfect agreement -with increasing approximation degree- with the analytical solution.  
%
%
Lastly, here we have added Table \ref{tab:Table_XHDG_Gaussian_temporal_p2_n32}, showing the pulse height at different time stamps; when the approximation functions used are second order i.e., $p=2$. Table \ref{tab:Table_XHDG_Gaussian_temporal_p2_n32} proves that XHDG discretization is accurate approximating the pulse height, independent of the flux definition used at various time stamps during the overall simulation. 

\begin{table}
\begin{center}
\caption{Temporal variation of the Gaussian pulse height evolving in XHDG discretized circular cut domain}\label{tab:Table_XHDG_Gaussian_temporal_p2_n32}
\vspace{0.3cm}
\resizebox{0.41\textwidth}{!}{%
\begin{tabular}{ccccc}
\hline
\hline
\multicolumn{2}{c}{central flux}&&\multicolumn{2}{c}{upwind flux}\\ 
\hline
$t$&$h$ && $t$&$h$\\
\hline
$0$&$1$&&$0$&$1$\\
$0.1$&$0.6341$&&$0.1$&$0.6335$\\
$1$&$0.2054$&&$1$&$0.2056$\\
$1.25$&$0.1606$&&$1.25$&$0.1608$\\
\hline
\hline
\end{tabular}}
\end{center}
\end{table}
\newpage
\section{Conclusions and final remarks}\label{sec:conclusions}

We have presented novel XHDG formulation for the solution of steady and time dependent solution of convection diffusion equations over unfitted domains. With unfitted XHDG discretization; the background mesh is cheap and easy to create and the domain of interest is defined with a level set function such that it can cut through the background mesh arbitrarily. With this unfitted strategy, the advantageous properties of HDG method, such as local and built in stabilization, element by element fashion, reduced number of degrees of freedom, optimal and super convergence is kept while removing the mesh-fitting-to-domain-boundary restriction; saving substantial mesh adaptation costs. 

In XHDG discretization, for the elements not cut by the domain boundary, standard HDG discretization is followed, while for cut elements a novel modified weak form is derived. Keeping the original unknown structure, to calculate the integrals over the cut elements and faces, a modified quadrature rule is used. A new unknown variable at the cut boundary is first introduced and then eliminated from the weak form imposing the boundary conditions. Extension to time dependent problems are as well presented, using first order Backward Euler time discretization scheme.    

Being different than other unfitted HDG strategies, with XHDG; the unknown structure of HDG is kept, high order $p\leq4$ convergence is proven, unsteady problems are treated and no extensions from domain or limitations on extension distance needed to be discussed. 

The XHDG discretization is tested over three steady and one unsteady problem where the cut boundary is either circular or in rather complex peanut shape. The convergence results of the numerical examples prove that while removing the mesh adaptation costs, XHDG shows optimal and super convergence for flow problems that is either diffusion or convection dominated.

Currently, XHDG formulation is being adapted for the solution of non-linear convection diffusion and Navier Stokes equations; together with a stabilization strategy to handle the ill conditioning issue at bad cut elements.  

\section{Acknowledgements}
This work is supported by The Scientific and Technological Research Council of T\"urkiye (T\"UBITAK), Career Development Program (CAREER), project no:121M947.


\newpage
\bibliographystyle{siam}
\bibliography{referencesXHDG}

\end{document}